\RequirePackage{fix-cm}
\documentclass[smallextended]{svjour3}       
\smartqed  
\usepackage{setspace}

\usepackage{amssymb}
\usepackage{amsmath}
\usepackage{amsfonts}
\usepackage{latexsym}
\usepackage{graphicx}
\usepackage{epstopdf}
\usepackage{caption}
\usepackage{subcaption} 
\usepackage[numbers,sort]{natbib}

\usepackage{float} 

\usepackage{hyperref} 

\journalname{Journal of Mathematical Biology}
\begin{document}

\title{Mathematical modelling of calcium signalling taking into account mechanical effects
}
%


\titlerunning{Calcium-Kaouri}        

\author{K. Kaouri      \and
        S. J. Chapman \and P. K. Maini
}


\institute{K.Kaouri \at
Department of Electrical Engineering, Computer Engineering, and Informatics\\
              Cyprus University of Technology, Cyprus \\
              \email{katerina.kaouri@cut.ac.cy}           
           \and
           S. J. Chapman \at
              Oxford Centre for Industrial and Applied Mathematics, Mathematical Institute, University of Oxford, UK \\
\and
          P. K. Maini\at
              Wolfson Centre for Mathematical Biology, Mathematical Institute, University of Oxford, UK
}

\date{Received: date / Accepted: date}

\maketitle

\begin{abstract}
Most of the calcium in the body is stored in bone. The rest is stored elsewhere, and calcium signalling is one of the most important mechanisms of information propagation in the body.  Yet,  many questions remain open. In this work, we initially consider the mathematical model proposed in Atri et al. \cite{ atri1993single}. Omitting diffusion, the model is a system of two nonlinear ordinary differential equations (ODEs) for the calcium concentration, and the fraction of  $IP_3$  receptors that have not been inactivated by the calcium.  We analyse in detail the system as the $IP_3$ concentration,  the \textit{bifurcation parameter}, increases presenting some new insights. We analyse asymptotically the relaxation oscillations of the model by exploiting a separation of timescales. Furthermore, motivated by experimental evidence that cells release calcium when mechanically stimulated and that,  in turn, calcium release affects the mechanical behaviour of cells, we propose an extension of the Atri model to a 3D nonlinear ODE mechanochemical model, where the additional equation, derived consistently from a full viscoelastic \emph{ansatz}, models the evolution of cell/tissue dilatation. Furthermore, in the calcium dynamics equation we introduce a new ``stretch-activation" source term that induces calcium release and which involves a new bifurcation parameter, the ``strength" of the source. Varying the two bifurcation parameters, we analyse in detail the interplay of the mechanical and the chemical effects, and we find that as the strength of the mechanical stimulus is increased, the $IP_3$ parameter range for which oscillations emerge decreases, until oscillations eventually vanish at a critical value. Finally, we analyse the model when the calcium dynamics are assumed faster than the dynamics of the other two variables. 

\end{abstract}


\section{Introduction}
\label{sec:Introduction}
Calcium signalling is one of the most important mechanisms of information propagation in the body, playing an important role as a second messenger in several processes such as muscle contraction, cardiac electrophysiology, bursting oscillations and secretion, synaptic plasticity, and adaptation in photoreceptors \cite{berridge2000, brini2009}. Therefore, understanding the mechanisms by which a cell controls its calcium concentration is of paramount importance. Through the sensing mechanisms of cells, external environmental stimuli can be transformed into intracellular or intercellular calcium signals which usually take the form of oscillations and waves. In particular, it is widely believed that the information carried by the environmental stimuli is encoded by the amplitude and/or frequency of calcium oscillations \cite{dupont2011, dupont2009, thul2014}. Therefore, understanding how the underlying mechanisms give rise to oscillations is of great scientific and physiological importance. 

Despite the versatility of calcium signals, in all cells calcium is stored and released from intracellular stores, such as the endoplasmic reticulum (ER), or the sarcoplasmic reticulum (SR). The concentration of calcium in these stores is 2-3 times larger than the concentration of calcium in the cytosol, and there is also a steep concentration gradient from the outside of the cell to the cytosol. These steep concentration gradients give the cells the advantage of being able to raise their calcium concentration quickly by opening calcium channels, relying only on passive flow down a steep gradient.

The release of calcium from the internal stores can be distinguished according to whether the release is dominated by the ryanodine receptors on the SR (mainly cardiac and skeletal muscle cells) or by inositol (1,4,5)-trisphosphate ($IP_3$) receptors on the ER (mainly non-muscle cells)--see \cite{keener1998, dupont2016}. The ryanodine receptors are activated by calcium and the $IP_3$ receptors are activated by $IP_3$; both types of receptors are modulated by the cytoplasmic calcium concentration, with calcium both activating and inactivating calcium release. Therefore both types of receptors are involved in a nonlinear feedback mechanism, frequently called Calcium Induced Calcium Release (CICR). The autocatalytic release of calcium stops when calcium level is sufficiently high; the inactivation level depends on the type of cell. The CICR mechanism has been used by several scientists to describe oscillations and waves of calcium inside cells\cite{keener1998}. The cytoplasm is an excitable medium with respect to calcium release in much the same way as the nerve axon is an excitable medium with respect to electrical potential \cite{keener1998, murray2001, grasman2012}.

Calcium signals, aside from being ubiquitous, are also diverse; depending on the type of cell they can vary by more than eight orders of magnitude in space and more than six orders of magnitude in time \cite{thul2014}. This wide spectrum of timescales and spatial scales makes calcium a fascinating object of study both for experimentalists and for mathematical modellers, and a great deal of work has been devoted to understanding many aspects of its behaviour and function. Yet, there are still many unanswered questions.

One of the most important and yet poorly understood issues in calcium signalling is the coupling of chemical and mechanical processes, i.e. how the cells sense mechanical stimuli and transform them into calcium signals and how calcium release affects the cells' mechanical properties and behaviour; this will be the major focus of this paper. Understanding these mechanochemical processes is important for understanding a wide range of processes, such as embryogenesis, wound healing, cancer modelling, mechanotransduction, and cell growth. 

It is well-known from experiments that calcium can be released due to local mechanical stimulation; mechanically stimulated calcium waves have been observed propagating through ciliated tracheal epithelial cells \cite{sanderson1981, sanderson1988, sanderson1990}, rat brain glial cells \cite{charles1991, charles1992, charles1993}, keratinocytes in the epidermis \cite{tsutsumi2009}, and many other cell types  \cite{young1999,beraheiter2005, yang2009,tsutsumi2009}. Also, in experiments with large tumours \cite{basson2015} it was observed that large extracellular pressure stimulated tumour proliferation via a mechanosensitive calcium channel. Moreover, spontaneous intercellular calcium waves have been observed in tissue slice preparations from the hippocampus \cite{dani1992} and in the intact liver \cite{robb1995}. Furthermore, in \cite{wallingford2001calcium} Wallingford et al. reported dramatic intercellular calcium waves in cells undergoing convergent extension in explants of gastrulating Xenopus embryos, which were often accompanied by a wave of contraction within the tissue.

Mathematical modelling has greatly advanced the understanding of calcium signalling and a multitude of deterministic and stochastic mathematical models have been developed by now -- see, for example,  
\cite{atri1993single,wilkins1998,goldberg2010,gracheva2001stochastic,sneyd1994model,sneyd1998calcium,timofeeva2003}; the comprehensive reviews \cite{sneyd2003, thul2014, rudiger2014stochastic} and the books \cite{keener1998, dupont2016}, among others. Due to the complexity of calcium signalling all models introduce approximations, depending on which aspect of calcium signalling they aim to examine. In this work we neglect stochastic effects. Even though the opening and closing of $IP_3$ receptors is a stochastic process, the deterministic models still have a good predictive power  in many cases, whilst being more amenable to analytical calculations \cite{thul2014, cao2014}. For a thorough discussion on deterministic vs. stochastic models see \cite{rudiger2014stochastic} and \cite{dupont2016}. 

Even though there are many models of calcium signalling, these behave in essentially similar ways, in that they model the CICR process in some way and give rise to oscillations and waves. (Tang et. al. have also shown in \cite{tang1996} that a number of seemingly different calcium models all have the same basic structure.) In this paper we start off with the model of Atri et al. \cite{atri1993single}, a two-dimensional system of ordinary differential equations (ODEs), and in section \ref{sec:Model} we present some new results on its bifurcation structure. We also study the model's asymptotic structure, under the assumption of a slow and a fast timescale at play, using Tikhonov's singular perturbation theory \cite{tikhonov1952}, frequently called geometric singular perturbation theory (GSPT) \cite{jones1995}. Understanding the asymptotic structure of ODE systems is important for practical applications, as it allows development of adequate approximate procedures that save computation time, and give intuitive insights, in a way that is not possible with exact models. GSPT reduces an ODE system of $k_1+k_2$ dimensions to lower-dimensional systems of dimensions $k_1$ and $k_2$ respectively, where the $k_1$-dimensional system is the \emph{slow} system and the  $k_2$-dimensional system is the \emph{fast} system \cite{tikhonov1952}. There are several applications of GSPT in the literature; for example, Suckley and collaborators \cite{suckley2003a, suckley2003b} have compared the asymptotic analysis of various excitable systems (Hodgkin-Huxley model, cardiac model) and Harvey et. al. \cite{harvey2011} used GSPT to analyse various intracellular calcium dynamics models. For the Atri model, when relaxation oscillations arise, we will see that the slow system is one nonlinear ODE, the fast system another nonlinear ODE; additionally, when $c$ becomes large the system reduces a  two-dimensional, linear system wvalid in a so-called `transition layer' \cite{grasman2012}. The presence of this transition layer makes the Atri model a nice, simple example of an excitable system which is qualitatively different to the classic Van der Pol oscillator and other similar equations such as the Fitzhugh-Nagumo system \cite{keener1998, grasman2012}.

In this work our main focus is mechanochemical modelling of calcium signalling, and in section \ref{sec:Mechanochemical} we introduce a new, mechanochemical, three-dimensional ODE model that combines the Atri model \cite{atri1993single} with a full viscoelastic \emph{ansatz} for modelling cell mechanics. To our knowledge the first mechanochemical models for biological problems were developed by Murray and collaborators in order to address questions in developmental biology-- see \cite{oster1989, murray1984, murray1988, murray2001}. The calcium dynamics assumed in the mechanochemical models of \cite{murray2001} are fairly basic and less sophisticated than the Atri model we consider here due to the lack of experimental data at that time; they assumed a bistable reaction diffusion model in which the application of stress can switch the calcium state from low to high stable concentration. Using the same bistable dynamics, Peradzynski in 2010 has analysed the latter models in two and in three dimensions, treating viscosity as a small parameter \cite{peradzynski2010}. Also, in 2010, Warren et al. constructed a complicated model which examines calcium wave propagation in mammalian airway epithelium when one of the cells is mechanically stimulated \cite{warren2010}. The model involves a reaction-diffusion PDE for calcium using the calcium dynamics of \cite{sneyd2004}, a reaction-diffusion PDE for $IP_3$ and a reaction-diffusion PDE for ATP, an enzyme which is released due to the mechanical stimulation of cells; they solved the resulting system numerically in 2D. In 2014 Kobayashi and collaborators modelled calcium waves in keratinocytes subjected to a local, instantaneous mechanical stimulus \cite{kobayashi2014}. Using the Atri model as a starting point, they developed an ODE-PDE model for the evolution of calcium, ATP, and $IP_3$ in each cell in a two-dimensional space. Although their model is simpler than the model in \cite{warren2010} they could only solve it numerically; they achieved satisfactory agreement with experimental results. In 2016 Yao and collaborators proposed a mathematical model where a network of mast cells (white blood cells) is in contact with interstitial flow \cite{yao2016}; they assumed one spatial dimension, and again a local, instantaneous mechanical stimulus. They derived a system of 4 ODEs and 1 PDE which they could only solve numerically; they also obtained partial agreement with the experiments in \cite{yang2009}. Summarising, the recent mechanochemical models in  \cite{warren2010, kobayashi2014, yao2016}, despite using more modern calcium dynamics than earlier models in \cite{oster1989, murray1984, murray1988, murray2001}, can only be analysed numerically. This severely limits our understanding of the underlying basic mechanochemical processes governing the interplay of calcium dynamics with cell mechanics. In constrast our proposed three-dimensional mechanochemical ODE model can be analysed with semi-analytical and asymptotic methods, and provides valuable new insights in the mechanochemical processes.

This paper is organised as follows. In section \ref{sec:Model} we present the Atri model  \cite{atri1993single} and analyse its rich bifurcation structure as the bifurcation parameter ($IP_3$ concentration) increases. Then, for the parameter range in which oscillations are observed we develop a new asymptotic analysis based on the GSPT, which has some differences from the asymptotic analysis of Van der Pol-like systems. In section \ref{sec:Mechanochemical} we present the new mechanochemical model which is a three-dimensional ODE system. To obtain the model we modify the Atri model by adding a new source term in the calcium dynamics equation modelling a `mechanical' stimulus, a source term that induces calcium release. The latter term involves a second bifurcation parameter, the strength of the mechanical stimulus. The third ODE models the evolution of a cell/tissue dilatation with the cell/tissue assumed to be a linear, Kelvin-Voigt material; the ODE is derived consistently from the full viscoelastic \emph{ansatz}. We then study the behaviour of the mechanochemical model in detail, and especially the oscillations arising.  Finally, for the parameter range for which the system has oscillations, assuming that calcium is a fast variable compared to the other two variables in the system, we develop a suitable asymptotic analysis. In section \ref{sec:Conclusions} we summarise the major points of our work and discuss possible extensions.

\section{The Atri model}
\label{sec:Model}

There are many models of calcium signalling - they all differ in the details but behave in essentially similar ways, as discussed above. In this paper we start-off with the Atri model \cite{atri1993single}, and analyse it in detail presenting some new, interesting, results on its bifurcation structure. The Atri model has been developed for intracellular spiral calcium waves in Xenopus oocytes but it has been used for modelling calcium dynamics in other types of cells too; for example, for glial cells in \cite{wilkins1998}, for mammalian spermatozoa in \cite{olson2010}, and for keratinocytes in \cite{kobayashi2014, kobayashi2016}. In \cite{wilkins1998} the Atri model has been used to study \emph{intercellular} calcium waves in glial cells. Modified Atri models have been developed in \cite{shi2008, harvey2011, liu2016}.  Therefore, the Atri model is a quite frequently used model of calcium dynamics in various cell types.

The Atri model equations are:
\begin{align}
\label{eq:Atria}
\frac{dc}{dt}&=J_{\rm flux}-J_{\rm pump}+J_{\rm leak}\\
\label{eq:Atrib}
\tau_h\frac{dh}{dt}&=\frac{k^2_2}{k^2_2+c^2}-h,\\
\nonumber
\textrm{where }J_{\rm flux}&=k_f\mu(p)h\frac{bk_1+c}{k_1+c},\,\,\,J_{\rm pump}=\frac{\gamma c}{k_{\gamma}+c},\,\,\,J_{\rm leak}=\beta,
\end{align}
where $c$ is the cytosolic calcium concentration, and $h$ is the fraction of $IP_3$ receptors on the ER that have not been inactivated by calcium.
In accordance to experimental evidence, the messenger $IP_3$ allows calcium to be released to the cytosol from an $IP_3$-sensitive store, in this case the ER. Furthermore, calcium is pumped out of the cytosol to the ER and also calcium leaks into the cytosol from outside the cell. Equation \eqref{eq:Atria} models the above processes: $J_{\rm flux}$ models the flux of calcium from the ER into the cytosol through the $IP_3$ receptors, assuming that calcium activates the $IP_3$ receptors quickly but inactivates them on a slower timescale; $\mu(p)=p/(k_{\mu}+p)$ is the fraction of $IP_3$ receptors that have bound $IP_3$ and is an increasing function of $p$, the  $IP_3$ concentration; $\mu$ will be used as the bifurcation parameter of the system. The constant $k_f$ denotes the calcium flux when all $IP_3$ receptors are open and activated. $J_{\rm pump}$ models the calcium pumped out of the cytoplasm back to the ER or out through the plasma membrane, and $J_{\rm leak}=\beta$ models the calcium leaking into the cytosol from outside the cell.  

In equation \eqref{eq:Atrib}, at steady state $h=\frac{k^2_2}{k^2_2+c^2}$, a decreasing function of $c$. Inactivation by calcium is assumed to act on a slow timescale characterised by the time constant $\tau_h$.  A detailed derivation of the model is presented in \cite{atri1993single} and the reader is referred there for more details. In this paper, we study the model in order to form a detailed understanding of the dynamics, and we present some new results on its bifurcation structure. 

We nondimensionalise the ODEs \eqref{eq:Atria}--\eqref{eq:Atrib} using $c=k_1\bar c$ and $t=\tau_h \bar t$, and drop bars for notational convenience. We obtain
\begin{align}
\label{eq:Atri1}
\frac{dc}{dt}&=\mu hK_1\frac{b+c}{1+c}-\frac{\Gamma c}{K+c}=F(c,h),\\
\label{eq:Atri2}
\frac{dh}{dt}&=\frac{K_2^2}{K_2^2+c^2}-h=G(c,h),
\end{align}
where $K_1=k_f\tau_h/k_1$, $K_2=k_2/k_1$, $\Gamma=\gamma \tau_h/k_1$, and $K=k_{\gamma}/k_1$. 
We use the parameter values in \cite{atri1993single}; $k_1=k_2=0.7\mu M$, $b=0.111$, $\gamma=2\mu M/s$, $k_f=16.2\mu M/s$, $k_{\gamma}=0.1\mu M$, $k_\mu=0.7\mu M$, $\tau_h=2s$, and hence $K_2=1$,  $\Gamma=5.71429$, and $K=1/7$. We note that the non-dimensional  parameter $K_1=46.285714$ and therefore the ratio of the timescale of calcium to the timescale of $h$, $1/K_1=0.021605$, is small, that is $c$ is a ``fast'' variable compared to $h$.  We will exploit this separation of timescales later on in order to analyse the Atri system asymptotically. 
 
 \subsection{Linear Stability}
\begin{figure}
\begin{center}
  \includegraphics[width=0.65\textwidth]{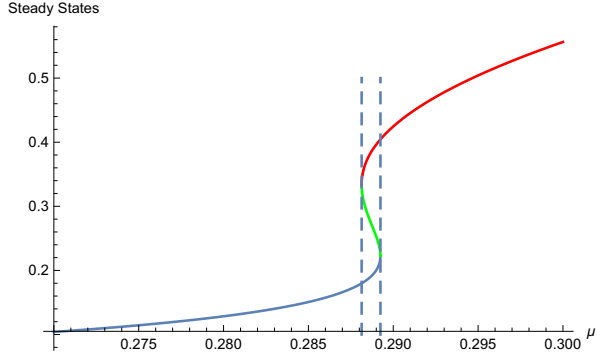}
    \caption{The steady states of \eqref{eq:Atri1}-\eqref{eq:Atri2} as a function of the bifurcation parameter $\mu$, as obtained from \eqref{eq:mu-SS}. As $\mu$ increases, from small to large there is one steady state, a double (degenerate) steady state at $\mu_1=0.28814$, then three steady states, then a double (degenerate) steady state at $\mu_2=0.28925$, and for values of $\mu$ larger than $0.28925$ one steady state.}
    \label{fig:pcFunctionmu}
\end{center} 
\end{figure}

The steady states (S.S.) of \eqref{eq:Atri1}-\eqref{eq:Atri2} are the intersections of the nullclines of the system. Setting
\begin{align}
\label{eq:Nullcline1}
F=0 \implies h&=\frac{\Gamma}{\mu K_1}\frac{c(1+c)}{(K+c)(b+c)},\\
\label{eq:Nullcline2}
G=0 \implies h&=\frac{1}{1+c^2}.
\end{align}
we obtain
\begin{align}
\label{eq:mu-SS}
\mu K_1\frac{1}{1+c^2}\frac{b+c}{1+c}-\frac{\Gamma c}{K+c}=0,
\end{align}
which can be cast as a quartic in $c$. In Figure  \ref{fig:pcFunctionmu} we plot the equilibrium curve \eqref{eq:mu-SS} in order to visualise the steady states' number and the corresponding value(s) of $c$, as $\mu$ is increased. 
The qualitative behaviour of the solutions of the system can be determined by plotting the nullclines \eqref{eq:Nullcline1} and  \eqref{eq:Nullcline2}. When the nullclines cross the system has a steady state, and when they touch the system has a double (degenerate) steady state. Nullcline \eqref{eq:Nullcline1} passes through the origin of the ($c$,$h$) plane, has a maximum at $h=h_M$ and saturates to the constant value $h=\frac{\Gamma}{\mu K_1}$ as $c \to \infty$. $h_M$ can be found analytically. Differentiating \eqref{eq:Nullcline1} we obtain
\begin{align}
\frac{dh}{dc}&=\frac{\Gamma}{\mu K_1}\frac{c^2(b+K-1)+2bcK+bK}{(K+c)^2(b+c)^2}
\label{eq:dhdc}
\end{align}
Setting $dh/dc=0$ leads to a quadratic in $c$. Re-arranging, and  since $c>0$, we discard the negative root, obtaining
\begin{align}
c_M&=\frac{bK+\sqrt{(1-b)b(K-K^2)}}{1-b-K}\label{eq:cM}\\
\intertext{and, hence, substituting \eqref{eq:cM} in \eqref{eq:Nullcline1} we obtain}
h_M&=h(c_M)=\frac{\Gamma}{\mu K_1}\frac{c_M(1+c_M)}{(K+c_M)(b+c_M)}.
\end{align}
Therefore, $h_M$ scales with $\Gamma/(\mu K_1)$ and depends on the parameters $K$, and $b$ in a more complicated manner. For the parameter values of \cite{atri1993single} we have $c_M=0.169=$constant and $h_M=0.279/\mu$.

Nullcline \eqref{eq:Nullcline2} is a decreasing function of $c$; it has a maximum at (0,1) and tends to $0$ as $c \to \infty$. For $\mu_1=0.28814$ and $\mu_2=0.28925$ the nullclines touch and there is a double steady state; for value of  $\mu<\mu_1$  and  $\mu>\mu_2$ there is one S.S. and there are three S.S. for $\mu_1<\mu <\mu_2$. The values of $\mu_1$ and $\mu_2$ are obtained by differentiating \eqref{eq:mu-SS} and finding the roots of $\frac{d\mu}{dc}=0$. Note: we present the values of $\mu$ with five decimal places because the bifurcation analysis depends sensitively on $\mu$, as we will see later.

\captionsetup[subfigure]{textfont=normalfont,singlelinecheck=on,justification=centering} 
We then linearise the system near the steady states. 
We determine the Trace (Tr), Determinant (Det) and Discriminant (Discr) of the Jacobian of the system as follows
\begin{align*}
&{\rm Tr}=F_{c}+G_h=F_c-1=\frac{\Gamma}{(K+c)}\left(-\frac{K}{K+c}+\frac{(1-b)\Gamma c}{(1+c)(b+c)}\right)-1\\
&{\rm Det}=F_cG_h-F_hG_c=-\frac{\Gamma}{(K+c)}\left(-\frac{K}{K+c}+\frac{(1-b)\Gamma c}{(1+c)(b+c)}\right)+\frac{2\Gamma c^2}{(1+c^2)(K+c)}\\
&{\rm Discr}=({\rm Tr})^2-{\rm 4Det}.
\end{align*}
We identify the bifurcations of the system as $\mu$ increases by determining where the Tr, Det, and Discr change sign. We easily do this by plotting Tr, Det, and Discr as functions of $c$ and using the numerical root-finding ``NSolve'' command in Mathematica. Below, we summarise the rich bifucation structure of the model as $\mu$ increases. This behaviour was not analysed in such detail in \cite{atri1993single} or later on in the literature, to the best of our knowledge.
\begin{itemize} 
\item[$\bullet$] $0<\mu<0.27828$: one stable node.
\item[$\bullet$] $\mu=0.27828$: the stable node becomes a stable spiral (bifurcation Discr=0)
\item[$\bullet$] $\mu=0.28814$: Stable spiral present. Also, a saddle and an unstable node (UN) emerge (bifurcation Det=0, fold point)
\item[$\bullet$] $\mu=0.28900$: the stable spiral becomes an unstable spiral. The other two S.S. are still a saddle and an unstable node. (Tr=0, Hopf bifurcation) 
\item[$\bullet$]  $\mu=0.28924$ the unstable spiral becomes an unstable node, and we have two unstable nodes and a saddle (Discr=0)
\item[$\bullet$] $\mu=0.28925$: one unstable node (Det=0, fold point)
\item[$\bullet$] $\mu=0.28950$: the unstable node becomes an unstable spiral (Discr=0)
\item[$\bullet$] $\mu=0.49500$: the unstable spiral becomes a stable spiral. (Tr=0, Hopf bifurcation)
\end{itemize}
From the various regimes above we are particularly interested in the ones corresponding to oscillations, due to their practical importance in calcium signalling. At $\mu=0.28900$ a Hopf bifurcation (HB) arises, the stable spiral becomes unstable, and we expect limit cycles in the nonlinear system. We also expect that these limit cycles will be sustained for a range of $\mu$ and then vanish close to $\mu=0.495000$ where the unstable spiral becomes a stable spiral.

\subsection{Numerical Solutions and Bifurcation Diagram}
As discussed in the Introduction (section \ref{sec:Introduction}) the amplitude and frequency of calcium oscillations encode signalling information so we are particularly interested in studying them further.   

\begin{figure}
\begin{subfigure}[t]{0.8\textwidth}
\includegraphics[width=1\linewidth]{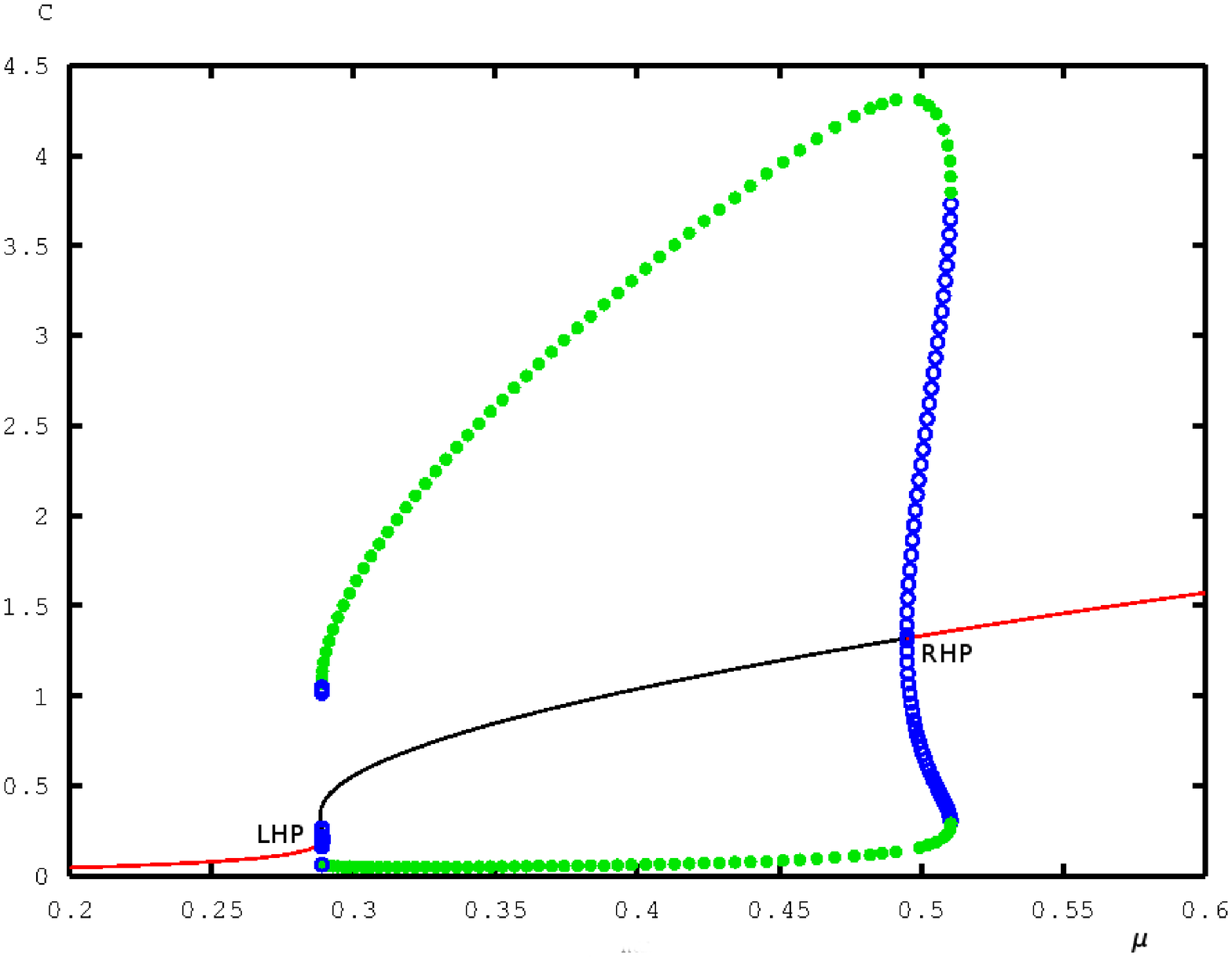}
\caption{}
\label{fig:Atri}
\end{subfigure}
\quad
\begin{subfigure}[t]{0.8\textwidth}
\includegraphics[width=1\linewidth]{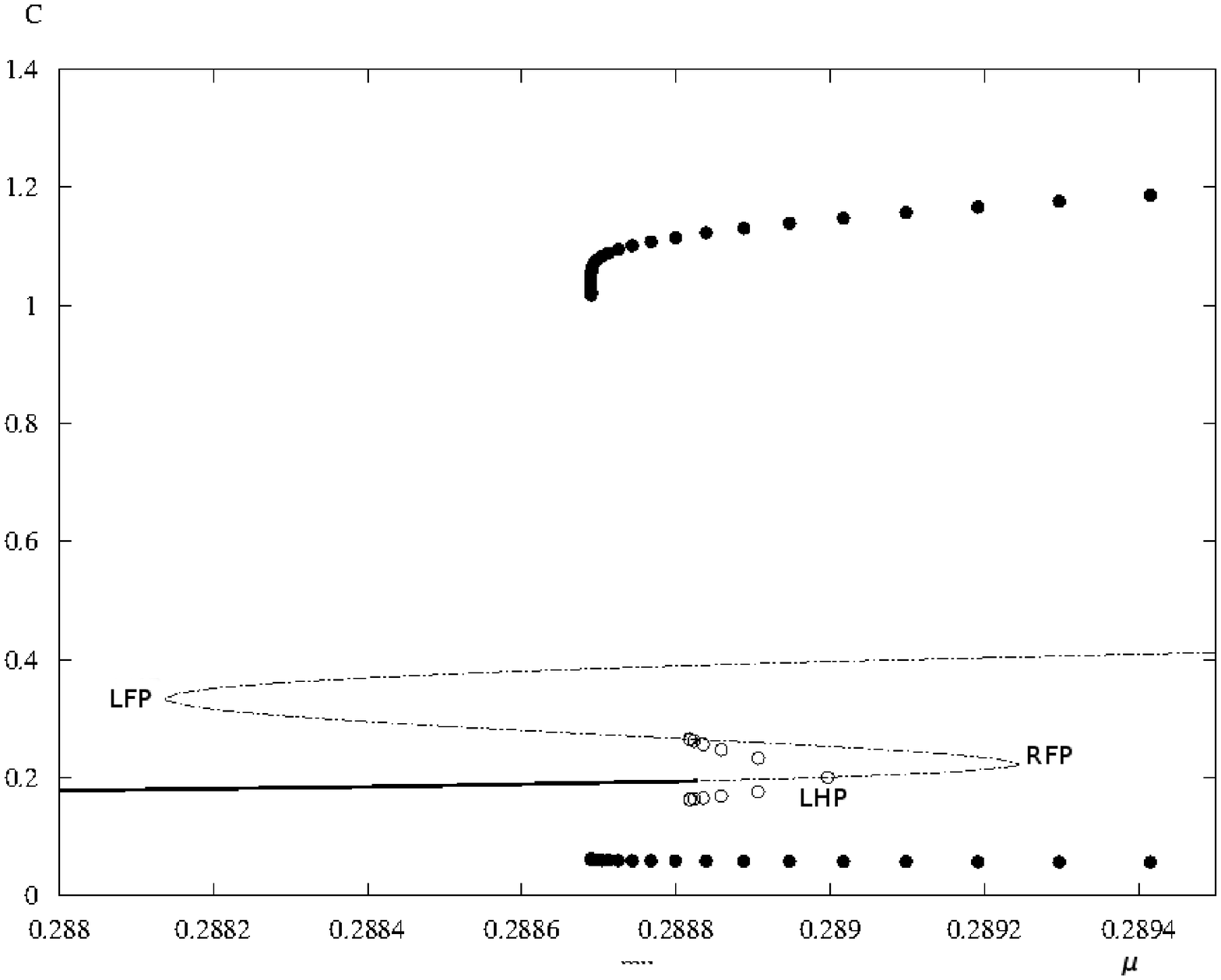}
\caption{}
\label{fig:AtriZoom}
\end{subfigure}
\quad
\begin{subfigure}[t]{0.8\textwidth}
\includegraphics[width=1\linewidth]{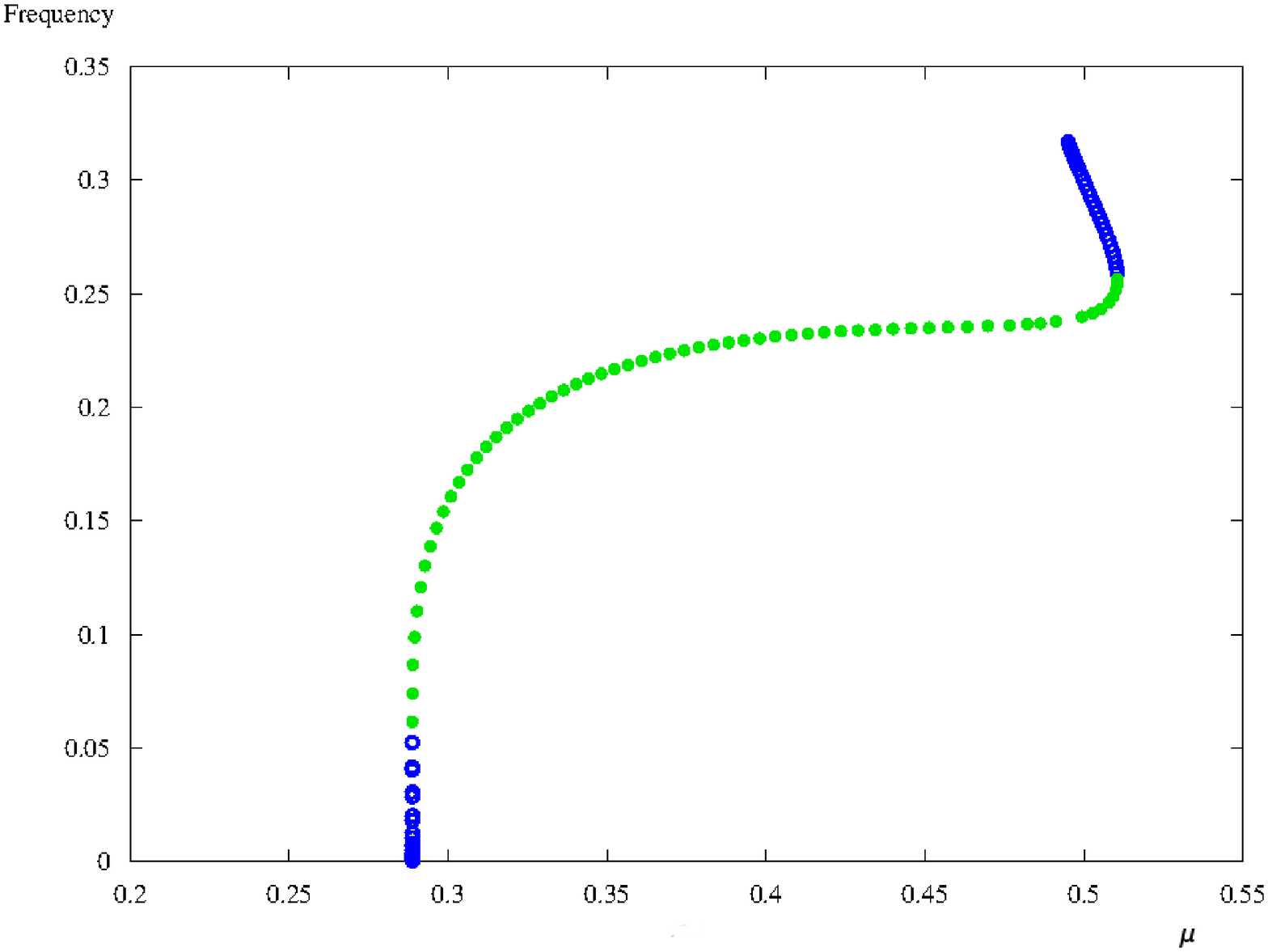}
\caption{}
\label{fig:AtriFrequency} 
\end{subfigure}
\caption{(a) Bifurcation diagram for the ODEs \eqref{eq:Atri1}--\eqref{eq:Atri2}, as $\mu$ changes - the solid circles represent stable limit cycles and the open circles represent unstable limit cycles--respectively with green and blue colour on the online version of the article). The LHP and the RHP are indicated. (b) Zooming in on the bifurcation diagram near the LHP. The LHP is very close to the RFP and a short, unstable, branch of periodic orbits emerges from it and terminates on the middle branch of the equilibrium curve. The stable limit cycle branch terminates to the right of the LFP.
(c) the frequency of the limit cycles.}
\label{fig:AtriBifnFrequency}
\end{figure}

In Figure \ref{fig:Atri}  we present the bifurcation diagram of the system. (We use AUTO as implemented in the continuation software XPPAUT by Ermentrout \cite{ermentrout2002simulating}). The two Hopf points at $\mu=0.28900$ and $\mu=0.49500$ signify the onset and suppression of oscillations respectively. There are stable limit cycles (represented by solid circles) and unstable limit cycles (represented by open circles).The branch of periodic orbits that comes out of the right Hopf point (RHP) at $\mu=0.49500$ is unstable and turns around at $\mu=0.5106$ where a stable and an unstable limit cycle coalesce (usually called the limit point of oscillations). (Note that only four decimal points are given by XPPAUT.) Therefore, for a small range of $\mu$ there is bistability between the stable limit cycle  and the stable upper steady state. The stable limit cycle branch undergoes a fold bifurcation where the frequency is zero and terminates to the right of the left fold point (LFP) in a homoclinic bifurcation  on the middle segment of the steady state curve. Also, from the left Hopf Point (LHP), which is very close to the right fold point (RFP), a very short branch of unstable periodic orbits emerges that terminates in a homoclinic bifurcation on the middle segment of the steady state curve. In Figure \ref{fig:AtriZoom} we zoom in on dynamics at the LHP. Note that near both Hopf points we see an extremely rapid growth in amplitude which makes computation of the exact bifurcation structure very sensitive to small changes in $\mu$, and also to the choice of the numerical method and stepsize used to compute the bifurcation diagram. In the intermediate region between the two Hopf points the amplitude of oscillations increases slowly with $\mu$. 

While the details of the dynamics in the very small interval of $\mu$ values near the LHP are mathematically interesting, from a biological point of view they are probably not so significant; in fact, for low $IP_3$ concentrations the validity of the deterministic model is questionable and a stochastic approach might be more appropriate. What is important here is the broad picture which the model has in common with other CICR models; for low and high $IP_3$ concentration there is an attracting steady state and for intermediate values there are stable oscillations (personal communication with James Sneyd and Vivien Kirk, and Chapter 5, p. 209 in the recent book by Dupont and collaborators \cite{dupont2016}).

In Figure \ref{fig:AtriFrequency} the frequencies of the stable and unstable limit cycles are shown with, respectively, solid and open circles; the range of $\mu$  for which both a stable and unstable limit cycle are sustained is clearly visible (double-valued part of the frequency curve). The limit point of oscillations at $\mu=0.5106$, where the stable and unstable limit cycle branches coalesce, is also clearly visible. The frequency of the stable limit cycles increases as $\mu$ increases and the lower (stable) branch approximates the square root of $\mu$ as predicted by bifurcation theory -- see \cite{kuznetsov2013} and other references on dynamical systems. This Atri model feature is compatible, qualitatively, with experiments in many cell types where the frequency of calcium oscillations increases with the $IP_3$ concentration --see, for example, \cite{berridge1988, berridge1989, parker1990}.

\begin{figure}[ht] 
  \begin{subfigure}[b]{0.5\linewidth}
    \centering
    \includegraphics[width=1\linewidth]{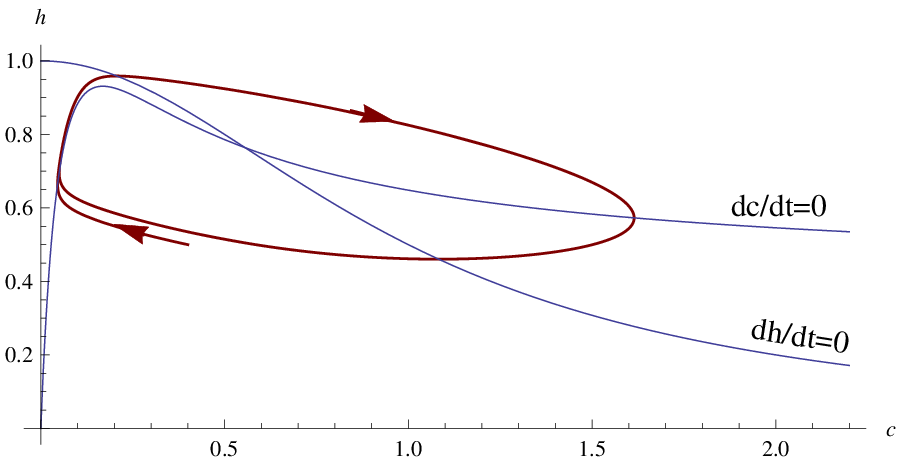}
    \caption{} 
  \label{fig:Solmu0Pt3a}
    \vspace{4ex}
  \end{subfigure}
  \begin{subfigure}[b]{0.5\linewidth}
    \centering
    \includegraphics[width=0.9\linewidth]{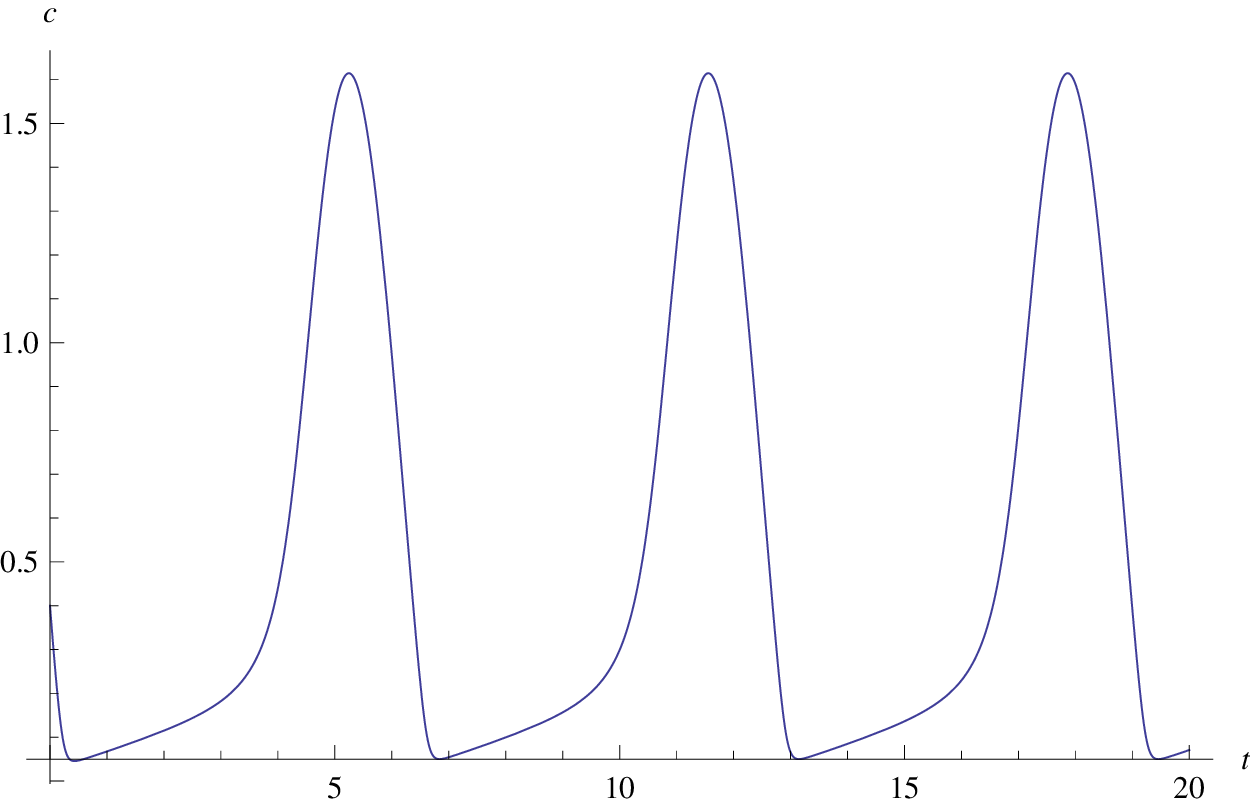}
    \caption{} 
   \label{fig:Solmu0Pt3b}
    \vspace{4ex}
  \end{subfigure} 
  \begin{subfigure}[b]{0.5\linewidth}
    \centering
    \includegraphics[width=1\linewidth]{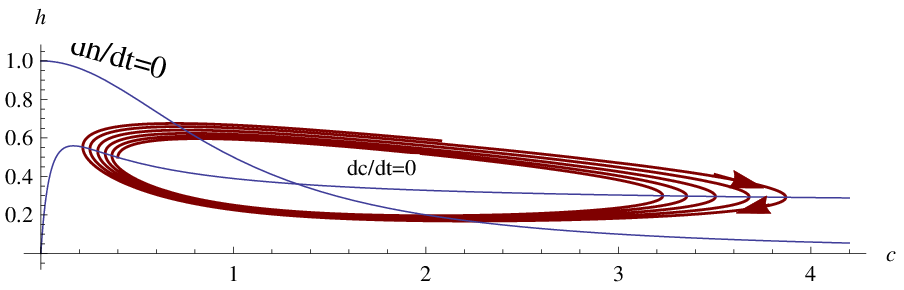}
    \caption{} 
  \label{fig:Solmu0Pt5a} 
  \end{subfigure}
  \begin{subfigure}[b]{0.5\linewidth}
    \centering
    \includegraphics[width=0.9\linewidth]{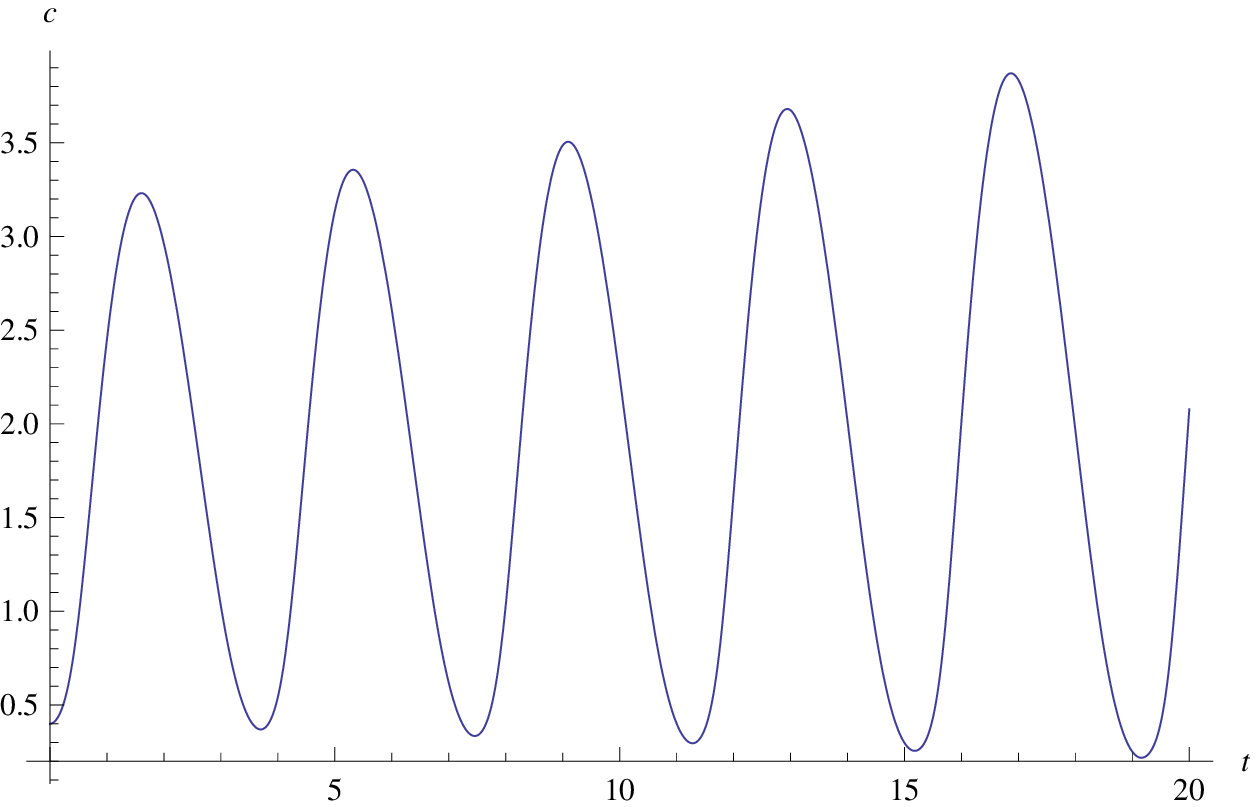}
    \caption{} 
   \label{fig:Solmu0Pt5b}
  \end{subfigure} 
  \begin{subfigure}[b]{0.5\linewidth}
    \centering
    \includegraphics[width=1\linewidth]{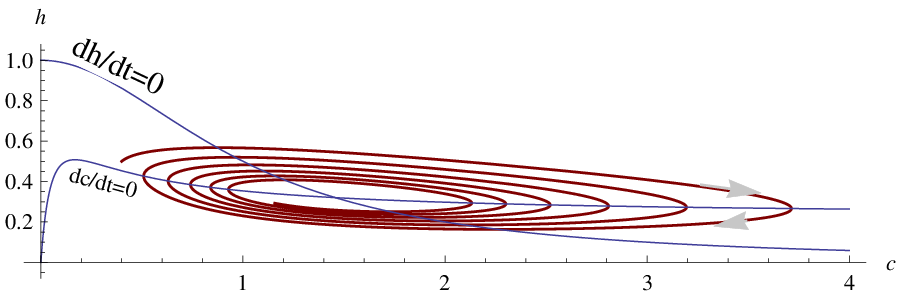}
    \caption{} 
  \label{fig:Solmu0Pt55a}
    \vspace{4ex}
  \end{subfigure}
  \begin{subfigure}[b]{0.5\linewidth}
    \centering
    \includegraphics[width=0.9\linewidth]{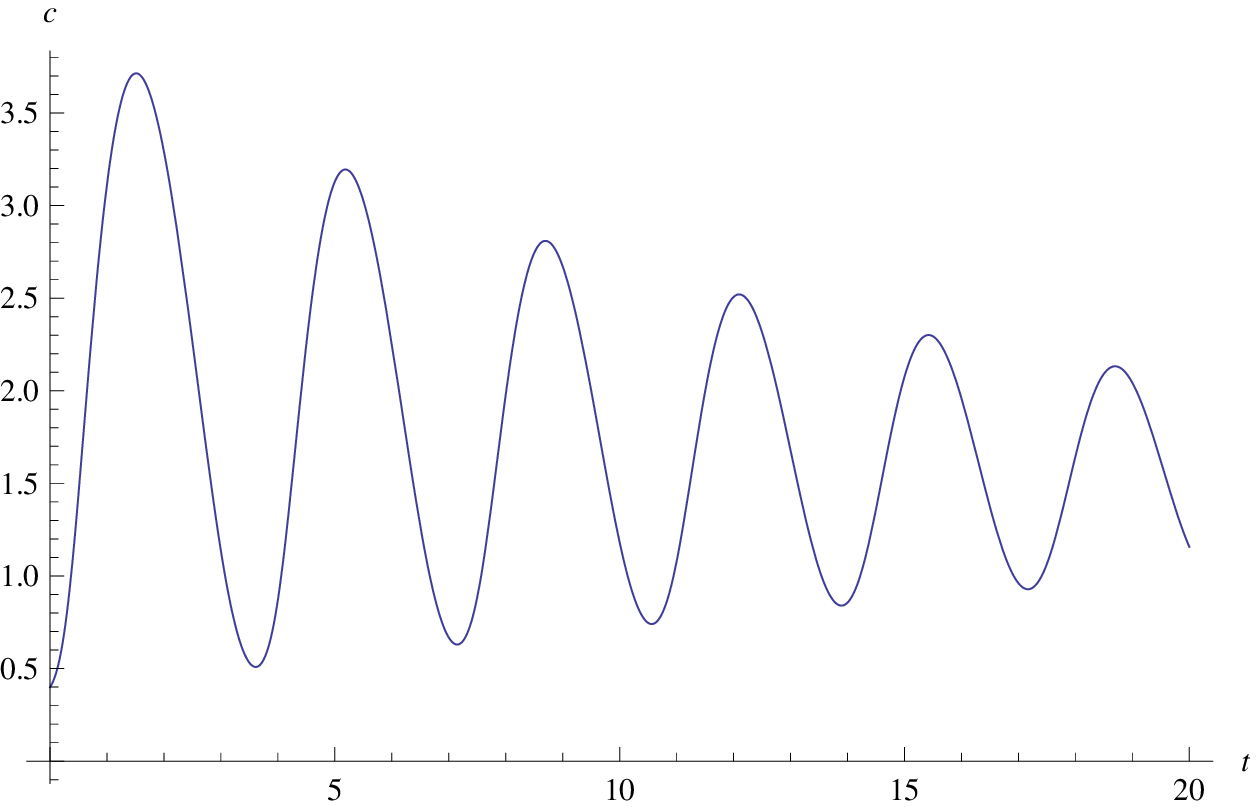}
    \caption{} 
   \label{fig:Solmu0Pt55b}
    \vspace{4ex}
  \end{subfigure} 
  \caption{Numerical solutions of ODEs \eqref{eq:Atri1}--\eqref{eq:Atri2} for initial conditions $c_0=0.4$, $h_0=0.5$ and nullclines for $\mu=0.3$ (a)-(b), $\mu=0.5$ (c)-(d), and $\mu=0.55$ (e)-(f). (Left column: Nullclines and the trajectory, Right column: evolution of calcium as a function of $t$, $c(t)$)}
\label{fig:Solmu}
\end{figure}

We then solve the nonlinear Atri ODEs \eqref{eq:Atri1}--\eqref{eq:Atri2} numerically and present the evolution of calcium with time for selected increasing values of $\mu$. Note that for  values of $\mu$ smaller than $\mu=0.28900$ the system exhibits excitable behaviour. For example, for $\mu=0.27$, with the initial condition $c_0=0.4$, $h_0=0.5$  the system settles quickly to the steady state but if we change the initial condition to $c_0=1$, $h_0=1$ the system exhibits an action potential,  i.e. it performs a large excursion before settling to equilibrium, a classic feature of \emph{excitable} systems \cite{keener1998, murray2001, grasman2012}.  For $\mu=0.3$, and with initial condition $c_0=0.4$, $h_0=0.5$, we see in Figures \ref{fig:Solmu0Pt3a}-\ref{fig:Solmu0Pt3b} that relaxation oscillations arise. For the same initial condition but increasing $\mu$ to $0.5$ relaxation oscillations persevere with increased amplitude and frequency (Figures \ref{fig:Solmu0Pt5a}-\ref{fig:Solmu0Pt5b}). Keeping $\mu=0.5$ and changing the initial condition to $c_0=0$, $h_0=0$ oscillations die out. This is because of the bistability between the stable and unstable limit cycles near the RHP. Finally, when $\mu$ is increased to $0.55$ oscillations die out for all initial conditions (see Figures \ref{fig:Solmu0Pt55a}-\ref{fig:Solmu0Pt55b}). 


\subsection{Asymptotic analysis of the Atri model}
The behaviour of the Atri model, as discussed earlier, is similar to the qualitative behaviour of other calcium models. However, in the Atri model the nullcline \eqref{eq:Nullcline1} saturates to the constant value $h=\Gamma/(\mu K_1)$ as $c \to \infty$ and \emph{it is not a cubic} as in Van derPol-like systems. This leads to some qualitatively different features in the asymptotic analysis of the Atri model, as we will see below.

In the Atri model, with the parameter values of \cite{atri1993single}, $c$ is a fast variable compared to $h$, as noted earlier. Here we exploit this separation of timescales to develop an asymptotic analysis, based on the geometric singular perturbation theory (GSPT). In usual perturbation theory notation, we let the small parameter, $\epsilon$, represent the ratio of the timescales of the dynamics of $c$ and $h$, and we rewrite \eqref{eq:Atri1} and \eqref{eq:Atri2} as follows: 
\begin{align}
\label{eq:AtriND1}
\epsilon \frac{dc}{dt}&=\mu h\frac{b+c}{1+c}-\frac{\Gamma}{K_1}\frac{c}{K+c}=\frac{F(c,h)}{K_1}=\epsilon F=f(c,h),\\
\label{eq:AtriND2}
\frac{dh}{dt}&=\frac{1}{1+c^2}-h=G(c,h).
\end{align}
Note that in Figure \ref{fig:Solmu}, $\epsilon$ is fixed to the value $1/K_1=0.021605$, according to the experimental values given in \cite{atri1993single},
 but for the purposes of the asymptotic analysis we are interested in the limit $\epsilon \to 0$, assuming $\frac{\Gamma}{K_1}$ is fixed and O(1). 

\subsubsection{Fast system}
The transformation of time $t=\epsilon \tau$ converts the system \eqref{eq:AtriND1}--\eqref{eq:AtriND2} into
\begin{align}
\label{eq:AtriND1Fast}
\frac{dc}{d\tau}&=\mu h\frac{b+c}{1+c}-\frac{\Gamma c}{K_1(K+c)},\\
\label{eq:AtriND2Fast}
\frac{dh}{d\tau}&=\epsilon\left(\frac{1}{1+c^2}-h\right).
\end{align}
Equations \eqref{eq:AtriND1}--\eqref{eq:AtriND2} are the slow system and \eqref{eq:AtriND1Fast}--\eqref{eq:AtriND2Fast} are the fast system. The slow system and the fast system are equivalent for finite $\epsilon>0$ but have different properties in the limit of $\epsilon \to 0$.

Setting $c=c_0(\tau)+O(\epsilon)$, $h=h_0(\tau)+O(\epsilon)$ in equations \eqref{eq:AtriND1Fast}-\eqref{eq:AtriND2Fast} the fast system, at leading order, is:
\begin{align}
\label{eq:AtriND1Fasta}
\frac{dc_0}{d\tau}&=\mu h_0\frac{b+c_0}{1+c_0}-\frac{1}{K_1}\frac{\Gamma c_0}{K+c_0},\\
\label{eq:AtriND2Fastb}
\frac{dh_0}{d\tau}&=0\implies h_0=h_{00}={\rm constant}.
\end{align}
Therefore, at leading order, from \eqref{eq:AtriND2Fastb} the slow variable $h$ is constant and only the fast variable $c$ varies according to equation \eqref{eq:AtriND1Fasta}. We call this regime \emph{Phase I}.
The non-intersecting straight lines $h_0=h_{00}={\rm constant}$ fill the whole of $R^2$ and constitute the \emph{fast foliation} of the space as they describe the evolution of the system on the fast timescale. Each of the lines constitutes a \emph{leaf} of the fast foliation.

\subsubsection{The slow manifold and slow dynamics}
To determine the leading order of the slow system we set $\epsilon=0$ in \eqref{eq:AtriND1}. We thus obtain the curve $f(c,h)=0$, the slow manifold (SM). The SM is separated into a stable and an unstable branch, represented respectively by the roots of the quadratic 
\begin{align}
\nonumber
&\mu hK_1(b+c)(K+c)-\Gamma c(1+c)=0\\
\nonumber
\Rightarrow &(\mu hK_1-\Gamma)c^2+(\mu hK_1(b+K)-\Gamma)c+\mu hK_1bK=0\\
\nonumber
\Rightarrow &c=c_{-}(h)\,\,\textrm{and}\,\,c=c_{+}(h),
\end{align}
where $c=c_{-}(h)$ is the smallest of the two roots.
The slow motion is obtained by solving the ODE \eqref{eq:AtriND2}, after substituting in it $c_{-}(h)$--this procedure is sometimes called an adiabatic elimination of the fast variable\footnote{Alternatively, we can obtain $c$ for the slow motion by changing variables in \eqref{eq:AtriND2} using $dh/dt=(dh/dc)dc/dt$, and solving the resulting ODE
\begin{align}
\label{eq:AtriND2Slow}
\frac{dc}{dt}=\frac{1}{\frac{dh}{dc}}\left(\frac{1}{1+c^2}-h(c) \right).
\end{align}
under the condition $dh/dc>0$; note that the RHS of \eqref{eq:AtriND2Slow} has a singularity when $dh/dc=0$.}.
When the trajectory reaches the maximum point of nullcline \eqref{eq:Nullcline1}, $(c_M, h_M)$, it flies off into a fast motion again. The maximum point is also frequently called a \emph{break point}. The condition $dh/dc=0$ is equivalent to the condition $\partial f/\partial c=0$. Also, the stable part of the SM satisfies $\partial f/\partial c<0$, which is equivalent to $dh/dc>0$, and the unstable part satisfies $\partial f/\partial c>0$, which is equivalent to $dh/dc<0$. Therefore, the SM is also the locus of the equilibria of the fast system. We call this slow motion regime \emph{Phase II}.

\subsubsection{Typical relaxation oscillation of the Atri system}
\begin{figure}
\captionsetup{justification=centering}
\begin{subfigure}[t]{1\textwidth}
\includegraphics[width=0.85\linewidth]{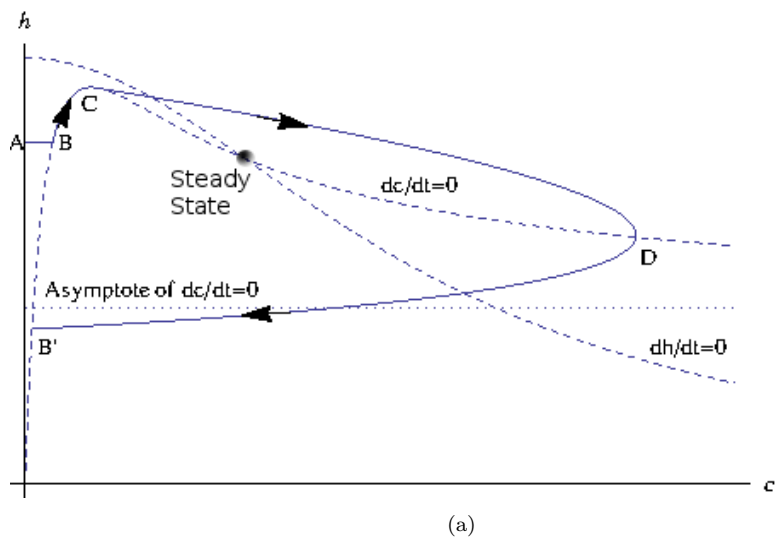}
\caption{}
\label{fig:AtriOscillationmu0Pt3}
\end{subfigure}
\quad
\begin{subfigure}[t]{1\textwidth}
\includegraphics[width=0.9\linewidth]{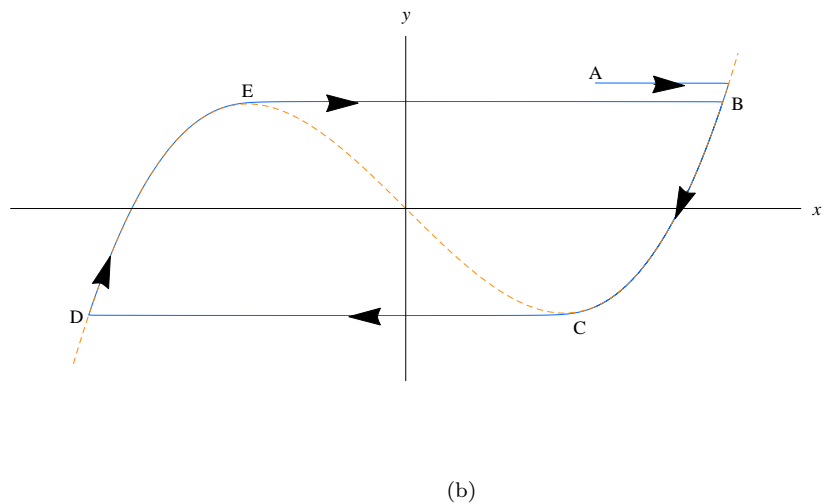}
\caption{}
\label{fig:VdPolOscillation} 
\end{subfigure}
\caption{(a) Typical relaxation oscillation for the Atri model \eqref{eq:AtriND1}--\eqref{eq:AtriND2}. The c-nullcline is not a cubic - it saturates to a constant value for $c$ large. (b) Typical relaxation oscillation for the Van der Pol oscillator which has a cubic nullcline.}
\label{fig:Atri-VdP}
\end{figure}

In Figure \ref{fig:AtriOscillationmu0Pt3} we plot a typical relaxation oscillation for the Atri system \eqref{eq:AtriND1}--\eqref{eq:AtriND2}, by piecing together the fast and slow parts. A fast motion takes place from A to B and a slow motion from B to C. 
The trajectory flies off at the maximum (break) point C and goes into a fast motion again, moving quickly towards D.  The motion from C to D, $c_0$ is governed by the reduced dynamics of equation \eqref{eq:AtriND1Fasta}, and since $\frac{dc_0}{d\tau}>0$ for $h_0>\frac{\Gamma}{\mu K_1}\frac{c_0(1+c_0)}{(K+c_0)(b+c_0)}$, $c_0$ is increasing when the trajectory is moving above the nullcline \eqref{eq:Nullcline1}. Since there is no other stable branch on the nullcline to which the trajectory can fly to, (as it would do in Var der Pol-like systems which have a ``cubic'' nullcline), $c$ inevitably becomes large, crosses the nullcline at D, and then starts moving towards the stable branch until it hits it again at the point B'. The region where $c$ gets large is called a `transition layer', according to the terminology of \cite{grasman2012}, and it necessitates the rescaling of $c$ with $\epsilon$. We shall call this regime \emph{Phase III} and analyse it further below.
After hitting the stable branch again at B' the trajectory goes up the SM in a slow motion again and proceeds to perform the same periodic motion indefinitely.
 
In Figure \ref{fig:VdPolOscillation}, in order to highlight the qualitative differences between the relaxation oscillations of the Atri model and of the Van der Pol oscillator we plot a typical relaxation oscillation for the latter\footnote{We use the equations $\frac{dx}{dt}=\frac{1}{\epsilon} (y+x-x^3/3), \frac{dy}{dt}=-\epsilon x$, with $\epsilon=0.025$.}.

\subsubsection{$c$ large region --`transition' layer}
To further investigate Phase III, we set $c=\hat c/\epsilon$, we assert the expansions $\hat c^{III}=\hat c_0^{III} (t)+O(\epsilon)$, $\hat h^{III}=\hat h_0^{III}(t)+O(\epsilon)$, and we substitute in the original system \eqref{eq:AtriND1}-\eqref{eq:AtriND2}. In the terminology of \cite{grasman2012} the latter region is called a `transition' layer. We find, dropping the superscript III,
\begin{align}
\label{eq:AtriND1-III}
\frac{d\hat c_0}{dt}&=\mu \hat h_0 - \frac{\Gamma}{K_1},\\
\label{eq:AtriND2-III}
\frac{d\hat h_0}{dt}&=-\hat h_0.
\end{align}
In the latter reduced system, which can be solved analytically, from \eqref{eq:AtriND1-III} $\displaystyle{\frac{d\hat c_0}{dt}=0}$ when $\hat h_0=\Gamma/\mu K_1$, $\displaystyle{\frac{d\hat c_0}{dt}>0}$ when $\hat h_0>\Gamma/\mu K_1$ and, correspondingly, $\displaystyle{\frac{d\hat c_0}{dt}<0}$ when $\hat h_0<\Gamma/\mu K_1$. This means that  the trajectory changes direction when it crosses the (asymptote of) the nullcline \eqref{eq:Nullcline1}, as expected (point D in Figure \ref{fig:AtriOscillationmu0Pt3})). After crossing the nullcline the trajectory starts moving towards the stable branch until it hits it again (point B' in Figure \ref{fig:AtriOscillationmu0Pt3}). Solving equations \eqref{eq:AtriND1-III}-\eqref{eq:AtriND2-III} we find
\begin{align}
\label{eq:AtriND1-III-Sol}
\hat c_0(t)&=-\mu \hat h_{00}e^{-t}-\frac{\Gamma}{K_1}t+d_0,\\
\label{eq:AtriND2-III-Sol}
\hat h_0(t)&=\hat h_{00}e^{-t}.
\end{align}
Matching Phase III with Phase II we find $\hat h_{00}=h_M$. To find the constant $d_0$ in \eqref{eq:AtriND1-III-Sol} we insert $\exp(-t)=1-t+...$ and match to Phase II, using $c \to 0$ as $t \to 0$. We find $d_0=\mu h_M$.
Therefore
\begin{align}
\label{eq:AtriND1IIISola}
\hat c_0(t)=\mu h_M (1-e^{-t})-\frac{\Gamma}{K_1}t.
\end{align}
To find the time needed for the trajectory to hit the nullcline for large $c$, $t_{\rm  TURNING}$, we set $\frac{d\hat c_0}{dt}=0$. We find 
\begin{align}
\label{eq:AtriND1IIISolb}
t=t_{\rm  TURNING}=ln\left(\frac{\mu K_1h_M}{\Gamma}\right).
\end{align}
Since $h_M$ scales with $\Gamma/(\mu K_1)$, for the parameter values we use $t_{\rm  TURNING}=\ln(2.2625)=0.82$, a constant.  Aside from $t_{\rm TURNING}$, other key feature of the relaxation oscillation are the peak of the oscillation, obtained by $c_{0MAX}=\hat c_0(t_{\rm TURNING})$, and $t_{\rm BACK}$, the time that the trajectory takes to hit the stable branch again, which is obtained by setting $\hat c_0=0$. From \eqref{eq:AtriND1IIISola} $t_{\rm BACK}$ is independent of $\mu$ also. Concluding, the dynamics in the transition layer are independent of $\mu$ at leading order.

Summarising, we have presented an asymptotic analysis of the Atri ODEs \eqref{eq:AtriND1}--\eqref{eq:AtriND2} when relaxation oscillations arise, as $\epsilon \to 0$, and we have plotted a typical relaxation oscillation which is qualitatively different from that of Van der Pol-like systems. At leading order, the 2D Atri system reduced to one nonlinear ODE on the slow manifold (Phase I), to one nonlinear ODE for the fast motion (Phase II), and to a simple 2D linear system in a novel, third region, called the `transition layer' (Phase III). The Atri model could be used as a nice textbook example for a nonlinear ODE system with a transition layer, in addition to the Van der Pol oscillator.

\section{A new mechanochemical model based on the Atri model}
\label{sec:Mechanochemical}
We now develop a new, simple, mechanochemical, ODE model which is suitable for modelling the mechanochemical processes in a single cell. The model can also be used to model mechanochemical processes in a thin, one-dimensional, monolayer of cells (called `tissue' from now on), that is assumed to be a continuum. Epithelial cells can spread or thicken when subjected to mechanical forces; through \emph{mechanosensing} (the process by which cells sense forces and transduce them into biochemical signals), the forces lead to calcium release from internal cell calcium stores \cite{lafaurie2016}. Also, the cells respond to an increase in intracellular free calcium by contracting. Therefore, there is a two-way mechanochemical feedback. The essential elements for building an appropriate mechanochemical model are a model for cell mechanics and a mechanism for calcium release, linked with a two-way feedback mechanism. The mechanisms governing cellular mechanotransduction are currently under intense investigation and outside the scope of this paper--more details can be found in the recent collection of Nature reviews \cite{naturereview2009}; another useful resource is the review by Barakat \cite{barakat2013}, examining mechanotransduction in epithelial cells involved in atherosclerosis. 


Our model consists of three ODEs; two of the ODES are the Atri model, where an additional, mechanical source term is added in the calcium equation. We also introduce a third ODE which governs the evolution of a new variable $\theta$, the cell/tissue dilatation, assuming that the cell/tissue is a  linear, one-dimensional, Kelvin-Voigt viscoelastic material. The model is 
\begin{align}
\label{eq:Atri3Da}
\frac{dc}{dt}&=\mu hK_1\frac{b+c}{1+c}-\frac{\Gamma c}{K+c}+\lambda \theta=R_1(c,\theta,h;\mu,\lambda),\\
\label{eq:Atri3Db}
\frac{d \theta}{dt}&=-\theta+T(c)=R_2(c,\theta),\\
\label{eq:Atri3Dc}
\frac{dh}{dt}&=\frac{1}{1+c^2}-h=R_3(c,h).
\end{align}
In \eqref{eq:Atri3Da} when $\theta$ increases, calcium is released, i.e. $\lambda \theta$ is a ``stretch-activation" source term, added to the original CICR term $\mu h K_1\frac{b+c}{1+c}$, and $\lambda\ge 0$ is a new parameter measuring the strength of this mechanical stimulus.  In equation \eqref{eq:Atri3Db} $T(c)\ge 0$ is a calcium-induced stress term that depends on the calcium concentration, and which causes an increase in $\theta$. 

Below we consistently derive \eqref{eq:Atri3Db} starting from a full viscoelastic \emph{ansatz} for the cell/tissue, as presented by Murray and collaborators in \cite{murray2001, oster1989, murray1984, murray1988}. A linear, viscoelastic, Kelvin-Voigt material with calcium-dependent cell traction/contraction is modelled as follows \cite{landau1986}:
\begin{align}
\label{eq:Atri3D-viscoelastic}
\nabla.\mathbf{\sigma}=0\Rightarrow \nabla.(\underbrace{\xi_1 \mathbf{e}_t+\xi_2 \theta_t \mathbf{I}}_\text{viscous stress}+\underbrace{E'(\mathbf{e}+\nu' \theta \mathbf{I})}_\text{elastic stress}-\underbrace{\tau(c)\mathbf{I})}_\text{cell traction/contraction}=0,
\end{align}
where $\mathbf{\sigma}$ is the stress tensor, $\mathbf{e}=\frac{1}{2}(\nabla \mathbf{u}+\nabla \mathbf{u}^T)$, with $\mathbf{u}$ the displacement vector,  is the strain tensor, $\theta=\nabla.\mathbf{u}$ is the dilatation, and $\mathbf{I}$ is the unit tensor; $\tau(c)$ is a calcium-dependent traction/contraction term. The constants $\xi_1, \xi_2$ are, respectively, the shear and bulk viscosities, and the constants $E'=E/(1+\nu)$ and $\nu'=\nu/(1-2\nu)$, where $E$ and $\nu$ are the Young's modulus and Poisson ratio, respectively. 
In one spatial dimension $\mathbf{e}=e=\theta=\frac{\partial u}{\partial x}$ and therefore \eqref{eq:Atri3D-viscoelastic} becomes, upon integrating with respect to $x$:
\begin{align}
\label{eq:Atri3D-viscoelastic-1D}
(\xi_1 +\xi_2) \theta_t+(1+\nu') \theta -\tau(c))=A,
\end{align}
where the constant of integration $A=0$ since $\tau=0$ when $c=0, \theta=0$ and $\theta_t=0$. Dividing \eqref{eq:Atri3D-viscoelastic-1D} by $\xi_1 +\xi_2$, rescaling according to
\begin{align}
\label{eq:Atri3D-viscoelastic-rescalings}
t^{\star}=\frac{1+\nu'}{(\xi_1 +\xi_2)}t, \,\,\,T(c)=\frac{\tau(c)}{1+\nu'},
\end{align}
and dropping the $\star$ from $t$, we deduce \eqref{eq:Atri3Db}.

In this work, we model $T(c)$ as an increasing function of $c$, as follows:
\begin{align}
\label{eq:Tmodel}
T(c)=\frac{\alpha c}{1+\alpha c},\,\,\,\,\,\,\,\,\alpha>0.
\end{align}
which satisfies the desired property that as $c=0$ the calcium-induced stress is $0$. As calcium levels increase we assume that the stress saturates to a constant value, $T_s=1$, i.e. the cell/tissue will be insensitive to further increases in calcium levels after a certain point. $T'(0)=\alpha$ is the rate of increase of $T$ at $c=0$ and $1/\alpha$ is the lengthscale of ascent to the saturation value $T_s$.  

\subsection{Linear stability analysis of the mechanochemical model when $\mu=0$}
Before proceeding to the analysis of the three-dimensional model we analyse the case $\mu=0$; biologically, this corresponds to switching off the $IP_3$ receptors and hence to no CICR flux. (Note that the $\lambda=0$ case is the Atri model.) The model \eqref{eq:Atri3Da}-\eqref{eq:Atri3Dc} then simplifies to
\begin{align}
\label{eq:Eqc-mu0}
\frac{dc}{dt}&=-\frac{\Gamma c}{K+c}+\lambda \theta,\\
\label{eq:Eqtheta-mu0}
\frac{d\theta}{dt}&=-\theta+T(c),\\
\label{eq:Eqh-mu0}
\frac{dh}{dt}&=\frac{1}{1+c^2}-h.
\end{align}
Equation \eqref{eq:Eqh-mu0} is decoupled from equations \eqref{eq:Eqc-mu0} and \eqref{eq:Eqtheta-mu0}, and we can thus continue with a two-dimensional analysis for equations \eqref{eq:Eqc-mu0}--\eqref{eq:Eqtheta-mu0}. The steady states satisfy $\displaystyle{\lambda T(c)=\frac{\Gamma c}{K+c}}$.
Evaluating the Jacobian of the 2D system \eqref{eq:Eqc-mu0}--\eqref{eq:Eqtheta-mu0}, we find that ${\rm Discr}>0, {\rm Tr}<0$ always, and that Det can be negative or positive\footnote{The entries of the Jacobian are
\begin{align*}
M_{11}=-\frac{\Gamma K}{(K+c^{\star})^2},\,\,M_{12}=\lambda,\,\,M_{21}=T'(c^{\star}),\,\,M_{22}=-1.
\end{align*}
Hence
\begin{align*}
{\rm Tr}=-\frac{\Gamma K}{(K+c^{\star})^2}-1<0,\,\,{\rm Det}=\frac{\Gamma K}{(K+c^{\star})^2}-\lambda T'(c^{\star}),\,\,{\rm Discr}=\left(\frac{\Gamma K}{(K+c^{\star})^2}-1 \right)^2+4\lambda T'(c^{\star})>0.
\end{align*}}. Therefore, for any $T(c)$, when $\mu=0$ the steady states can only be stable nodes or saddles and oscillations cannot be sustained. For some choices of $T$ a non-zero steady state may exist and this means biologically that even without an $IP_3$-induced calcium flux from the ER a non-zero calcium concentration can be sustained in the cytosol due to the stress-induced calcium release.

For $T(c)$ given by \eqref{eq:Tmodel}, there is always a steady state ($c^{\star},\theta^{\star}$)=($0,0$). A second steady state 
\begin{align}
\label{eq:SS2}
c^{\star}=\frac{\delta-K}{1-\alpha \delta},\,\,\textrm{ where } \delta=\frac{\Gamma}{\alpha \lambda},
\end{align}
exists only if
\begin{align}
\label{eq:SS2ConditionsA}
&\delta>K\,\, \textrm{and}\,\, \alpha \delta<1\,\Rightarrow \lambda<\frac{\Gamma}{\alpha K}\,\, \textrm{and}\,\, \lambda > \Gamma   \\
\label{eq:SS2ConditionsB}
\textrm{ or }\,\,&\delta< K\,\,\textrm{ and }\,\,\alpha \delta>1\,\Rightarrow \lambda>\frac{\Gamma}{\alpha K}\,\, \textrm{and}\,\, \lambda < \Gamma  .
\end{align}
The S.S. (0,0) loses its stability and becomes a saddle when the second stable S.S. emerges (stable node).  

\subsection{Linear stability analysis of the mechanochemical model}
\label{sec:LinStab}
The steady states of the system \eqref{eq:Atri3Da}--\eqref{eq:Atri3Dc} satisfy $R_1=0$, $R_2=0$, $R_3=0$, and hence satisfy into the following relationship:
\begin{align}
\label{eq:mu-SS-3D}
\mu K_1\frac{1}{1+c^2}\frac{b+c}{1+c}-\frac{\Gamma c}{K+c}+\lambda T(c)=0.
\end{align}
Using \eqref{eq:mu-SS-3D}, for any $T(c)$, we can easily plot the equilibrium surface, i.e. the steady states vs. $\mu$ and $\lambda$. We note that \eqref{eq:mu-SS-3D} reduces to \eqref{eq:mu-SS} for $\lambda=0$, as expected.

We now determine the stability of the steady states for any $T(c)$. The Jacobian of \eqref{eq:Atri3Da}--\eqref{eq:Atri3Dc}, at any S.S.,  is given by
\begin{align}
M=
  \begin{bmatrix}
    R_{1c} &  R_{1\theta} &R_{1h} \\
     R_{2c}& R_{2\theta}& R_{2h}  \\
    R_{3c}& R_{3\theta} & R_{3h} 
  \end{bmatrix}
=
  \begin{bmatrix}
    R_{1c} & \lambda& R_{1h}  \\
      T'(c)& -1 & 0 \\
   R_{3c}  & 0& -1
  \end{bmatrix}.
\end{align}
The characteristic polynomial is conveniently factorised 
\begin{align}
\label{eq:factorisedpoly}
(1+\omega)(\lambda T'(c)+(R_{1c}-\omega)(1+\omega)+R_{1h}R_{3c})=0,
\end{align}
where $\omega$ represents the eigenvalues. The factorisation is possible because: (i) \eqref{eq:Atri3Db},does not depend on $h$ explicitly and \eqref{eq:Atri3Dc} does not depend on $\theta$ explicitly, and hence $R_{2h}=R_{3\theta}$ 
(ii) $R_{2\theta}=R_{3h}$=-1. 

Solving \eqref{eq:factorisedpoly}, one eigenvalue is always equal to -1, and the linear stability, for any $T(c)$, can be studied through the \emph{quadratic}
\begin{align}
\omega^2-\omega(R_{1c}-1)-R_{1c}-R_{1h}R_{3c}-\lambda T'(c)=0.
\end{align}
We will use the latter quadratic to easily study the bifurcations of the system in the $\mu$-$\lambda$ parameter space, for any $T(c)$. We are especially interested in the amplitude and frequency of oscillations as $\lambda$ increases. Therefore,first, we are going to study the Hopf bifurcations of the system, which satisfy Tr=$0$, Det$>0$, Discr$<0$: 
\begin{align}
\label{eq:Tr}
\textrm{Tr}=0\implies R_{1c}-1=0\implies \frac{\mu K_1(1-b)}{(1+c^2)(1+c)^2}-\frac{\Gamma K}{(K+c)^2}-1=0.
\end{align}
(Note that when $\mu=0$, $\textrm{Tr}<0$ always.) Since \eqref{eq:Tr} and \eqref{eq:mu-SS-3D}
are linear in  $\mu$ and $\lambda$ we can easily express  $\mu$ and $\lambda$ in terms of $c$ as follows:
\begin{align}
\label{eq:MUvsC}
\mu(c)&=\frac{(1+c^2)(1+c)^2}{K_1(1-b)}\left(1+\frac{\Gamma K}{(K+c)^2}\right),\\
\label{eq:LAMBDAvsC}
\lambda(c)&=\frac{1}{T(c)}\left(\frac{\Gamma c}{K+c}-\frac{(b+c)(1+c)}{1-b}\left(1+\frac{\Gamma K}{(K+c)^2}\right)\right).
\end{align}
We observe that the function $T(c)$ appears only in \eqref{eq:LAMBDAvsC}.
We can, thus, easily plot \eqref{eq:MUvsC}, \eqref{eq:LAMBDAvsC} parametrically, and obtain the \emph{Hopf curve}, for \emph{any} $T(c)$. The interior of the Hopf curve corresponds to an unstable spiral and approximates the $\mu$-$\lambda$ region corresponding to relaxation oscillations (limit cycles) in the full nonlinear system. 

We can also determine analytically a parametric expression for the \emph{fold curve} in the $\mu$-$\lambda$ plane; setting
\begin{align}
\nonumber
&\textrm{Det}=0\implies -R_{1c}-R_{1h}R_{3c}-\lambda T'(c)=0\\
\label{eq:Det1}
& \implies \mu\frac{K_1}{(1+c)(1+c^2)}\left(\frac{1-b}{1+c}-2c(b+c) \right)+\lambda T'(c)=\frac{\Gamma K}{(K+c)^2}.
\end{align}
Equations \eqref{eq:Det1} and \eqref{eq:mu-SS-3D} constitute a linear system for $\mu$ and $\lambda$, and we can again easily derive parametric expressions for $\mu$ and $\lambda$ in terms of $c$. For $\lambda=0$ the results agree with the Atri model, as expected. 

Similarly, to determine parametric expressions for the curve on which $\rm Discr$ changes sign, we set
\begin{align}
\label{eq:Discr}
&{\rm Discr}=0\implies (R_{1c}+1)^2+4R_{1h}R_{3c}+4\lambda T'(c)=0,
\end{align}
which is quadratic in $\mu$ and linear in $\lambda$. Combining \eqref{eq:Discr} with \eqref{eq:mu-SS-3D} we can determine again parametric expressions for $\mu$ and $\lambda$.
In summary, we have developed a fast method for determining the key curves of the bifurcation diagram, for \emph{any} $T(c)$. Below we consider a simple, illustrative example.

\subsubsection{Two-parameter bifurcation diagram for $\displaystyle{T=\frac{\alpha c}{1+\alpha c}}$}
\begin{figure}
\begin{center}
\includegraphics[width=0.85\textwidth]{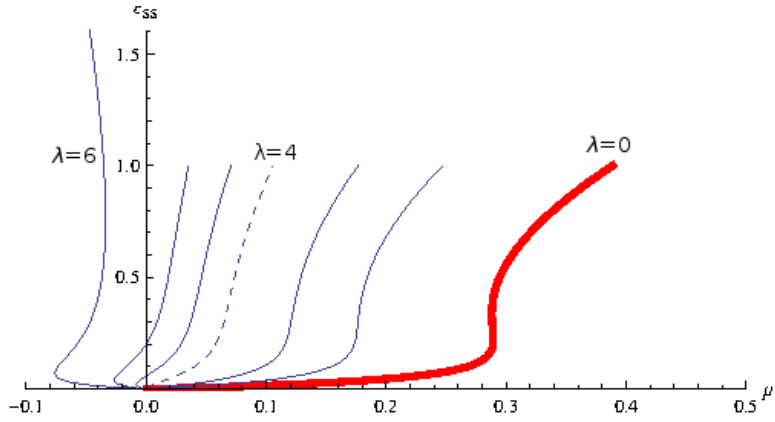}
\caption{Steady states of the system \eqref{eq:Atri3Da}--\eqref{eq:Atri3Dc} as $\mu$ is increased, for $\lambda=0, 2, 3, 4, 4.5, 5, 6$, when $T=10c/(1+10c)$. $\lambda$ increases from right to left -- the thick curve is 
$\lambda=0$ and the dashed curve is $\lambda=4$.}
\label{fig:pBifDiagfmuLam}
\end{center}
\end{figure}

We choose  $T=10c/(1+10c)$, and in Figure \ref{fig:pBifDiagfmuLam} we use \eqref{eq:mu-SS-3D} to plot the steady state vs. $\mu$ for selected, increasing values of $\lambda$. 
We see that for $\lambda<4$ the bifurcation diagram is qualitatively similar to the $\lambda=0$ diagram.
However, at $\lambda=4$ the qualitative behaviour of the equilibrium curve changes and a second non-zero steady state appears for $4<\lambda<40/7=5.71429$, as predicted by the condition \eqref{eq:SS2ConditionsB} when we set $\alpha=10$, $K=1/7$ and $\Gamma=5.71429$. Note that for the aforementioned range of $\lambda$ the part of the equilibrium curve corresponding to negative values of $\mu$ is not biologically relevant.

\begin{figure}	
	\centering
	\begin{subfigure}[t]{5.3in}
		\centering
		\includegraphics[width=5.0in]{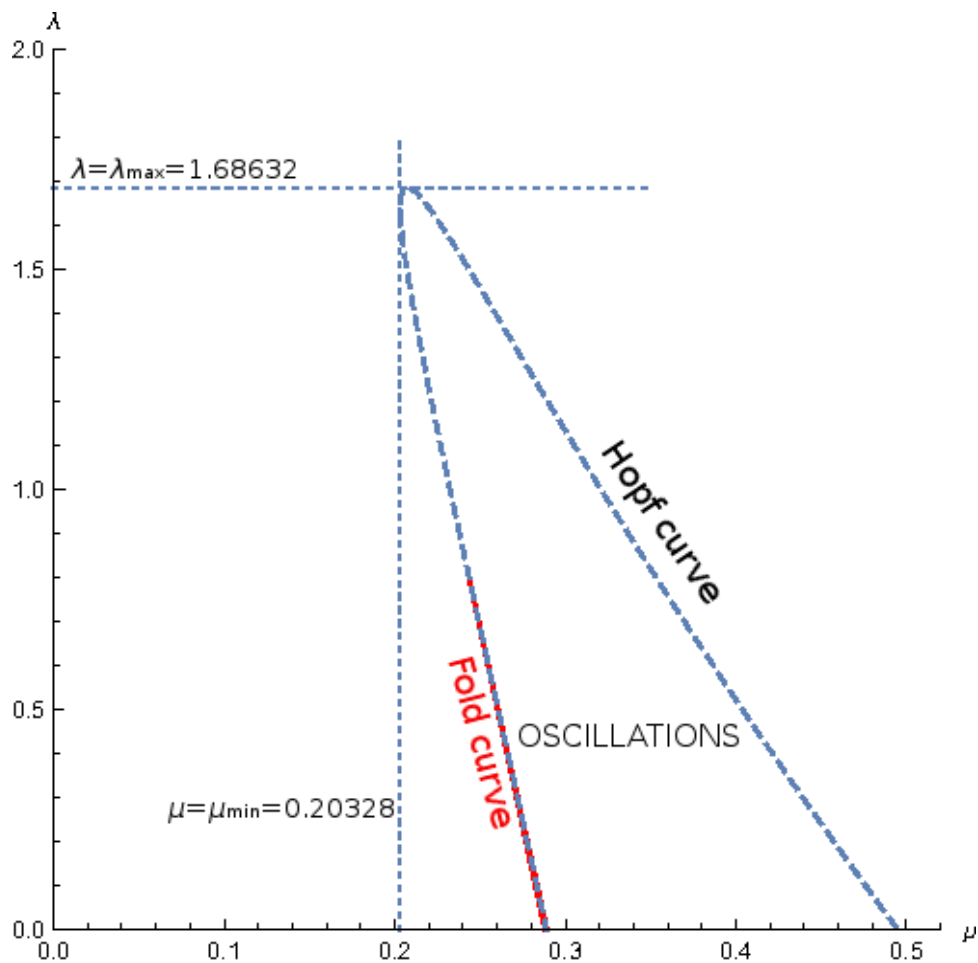}
		\caption{}\label{fig:HopfCurveFold}		
	\end{subfigure}
	\quad
	\begin{subfigure}[t]{5.3in}
		\centering
		\includegraphics[width=3.2in]{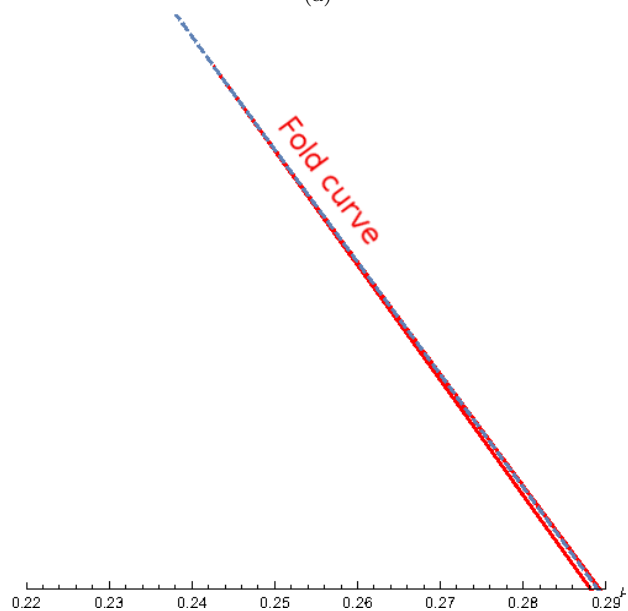}
		\caption{}\label{fig:HopfcurveFoldZoom}
	\end{subfigure}
	\caption{$\mu$-$\lambda$ bifurcation diagrams of the system \eqref{eq:Atri3Da}--\eqref{eq:Atri3Dc} when $T(c)=10c/(1+10c)$: a) the Hopf curve (dashed-blue colour in online version) and fold curve (solid-red colour in online version); the two branches of the fold curve are very close to each other. The horizontal and vertical dashed lines correspond, respectively, to the maximum value of $\lambda$, $\lambda_{\rm \max}=1.68632$ and to the minimum value of $\mu$, $\mu_{\rm min}=0.20328$ for which oscillations are sustained.   (b) Zooming in on the fold curve of Figure \ref{fig:HopfCurveFold}. Inside the fold curve there are \emph{three} steady states, there is a double steady state on the curve, and \emph{one} steady state outside the curve.}\label{fig:AtriMechsBDiag}
\end{figure}

In Figure \ref{fig:HopfCurveFold} we plot the Hopf curve and the fold curve for $T(c)=\frac{10c}{1+10c}$. We observe the following: (i) for $\lambda=0$ we recover the Hopf points of the Atri model. (ii) As $\lambda$ increases the  $\mu$--range of oscillations decreases. (iii) For a critical maximum value of $\lambda$, $\lambda_{\rm \max}$, the oscillations are suppressed. Solving  $\frac{d\lambda}{dc}=0$, using \eqref{eq:LAMBDAvsC}, we can find $\lambda_{\rm \max}=1.68632$, and we can use \eqref{eq:MUvsC} to determine $\mu_{\rm \max}=0.20735$. (iv) Also, from \eqref{eq:MUvsC}, setting $\frac{d\mu}{dc}=0$, and using \eqref{eq:LAMBDAvsC} we also find $\mu_{\rm min}=0.20328$, the minimum value of $\mu$ for which oscillations can be sustained ($\lambda_{\rm min}=1.63989$). The fold curve consists of two branches very close to each other; the interior of the fold curve corresponds to three steady states and the exterior to one steady state (so we have bistability for an increasingly smaller range of $\mu$ as $\lambda$ increases) until eventually at $\lambda \approx 0.83$ the two branches merge; for values of $\lambda \ge 0.83$ the system has only one steady state. In Figure \ref{fig:HopfcurveFoldZoom} we zoom in on the two branches of the fold curve.

\begin{figure}	
	\centering
	\begin{subfigure}[t]{4.6in}
		\centering
		\includegraphics[width=4.6in]{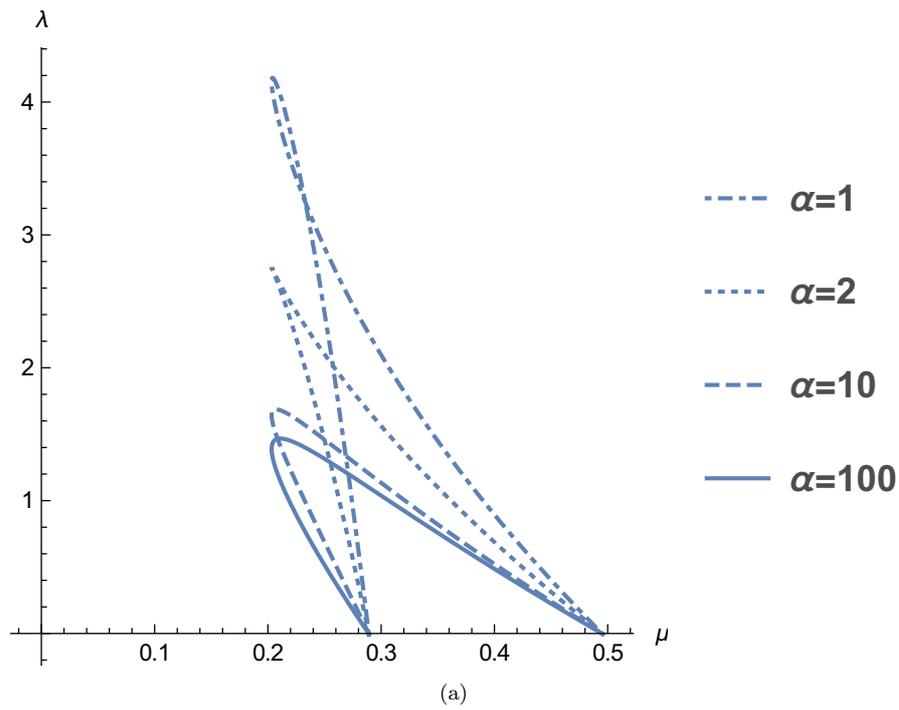}
		\caption{}\label{fig:HopfCurvesAlphaVaries}
	\end{subfigure}
	\quad
	\begin{subfigure}[t]{3.6in}
		\centering
		\includegraphics[width=3.6in]{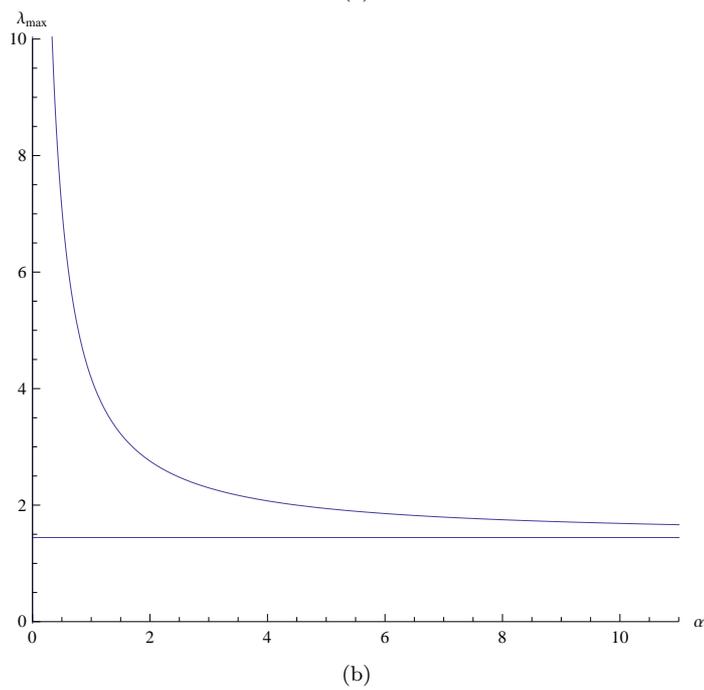}
		\caption{}\label{fig:pLamCritvsAlpha}
	\end{subfigure}
	\caption{(a) Hopf curves for the system \eqref{eq:Atri3Da}--\eqref{eq:Atri3Dc} and $T(c)=\frac{\alpha c}{1+\alpha c}$, $\alpha=1, 2, 10, 100$ (b) The maximum value of $\lambda$ for which oscillations are observed, $\lambda_{\rm max}$, as a function of $\alpha$.}\label{fig:}
\end{figure}

We are also interested in how/if the Hopf curve changes qualitatively as we vary the parameter $\alpha$, equal to the rate of growth of $T$ at $c=0$. We choose four different cases of $T(c)$ and in Figure \ref{fig:HopfCurvesAlphaVaries} we plot the Hopf curves; we see that they are qualitatively similar for $\alpha=10$ and $\alpha=100$ but that as $\alpha$ decreases they change qualitatively; for $\alpha \approx 2$ the Hopf curve develops a cusp, and for smaller values of $\alpha$ we see a ``bow-tie''. This is another bifurcation, with $\alpha$ as a bifurcation parameter. The cusped Hopf curve is the bifurcation curve of a cusp catastrophe surface, according to the catastrophe theory developed by Zeeman \cite{zeeman1977}, and subsequently by Stewart and collaborators \cite{poston2014, stewart2014}. In summary, we see that for any $\alpha$ oscillations vanish after a certain value of $\lambda$, denoted by $\lambda_{\rm \max}$; and as $\alpha$ is increased $\lambda_{\rm \max}$ decreases. Also, since \eqref{eq:MUvsC} does not depend on $T(c)$,  the minimum value of $\mu$ for which oscillations can be observed is equal to $\mu_{\rm min}=0.20328$, for any $\alpha$. 

To investigate more systematically the decrease of $\lambda_{\rm \max}$ with $\alpha$ we determine parametric expressions for $\lambda_{\rm \max}$ and $\alpha$ as functions of $c$, and we plot $\lambda_{\rm \max}$ vs. $\alpha$ in Figure \ref{fig:pLamCritvsAlpha}. We see in Figure  \ref{fig:pLamCritvsAlpha} that as $\alpha$ increases, $\lambda_{\rm max}$ decreases monotonically, and hence oscillations are sustained for a smaller range of $\lambda$ (as indicated in Figure \ref{fig:HopfCurvesAlphaVaries}). Also, since 
$\lambda_{\rm max}(\alpha)$ asymptotes to a positive value as $\alpha \to \infty$, our model, for any $T(c)=\frac{\alpha c}{1+\alpha c}$, predicts that there will always be oscillations for a certain range of $\mu$ values. 

\subsubsection{Hopf curves for $\displaystyle{T=\frac{\alpha c^2}{1+\alpha c^2}}$}
\begin{figure}
\begin{center}
\includegraphics[width=0.45\textwidth]{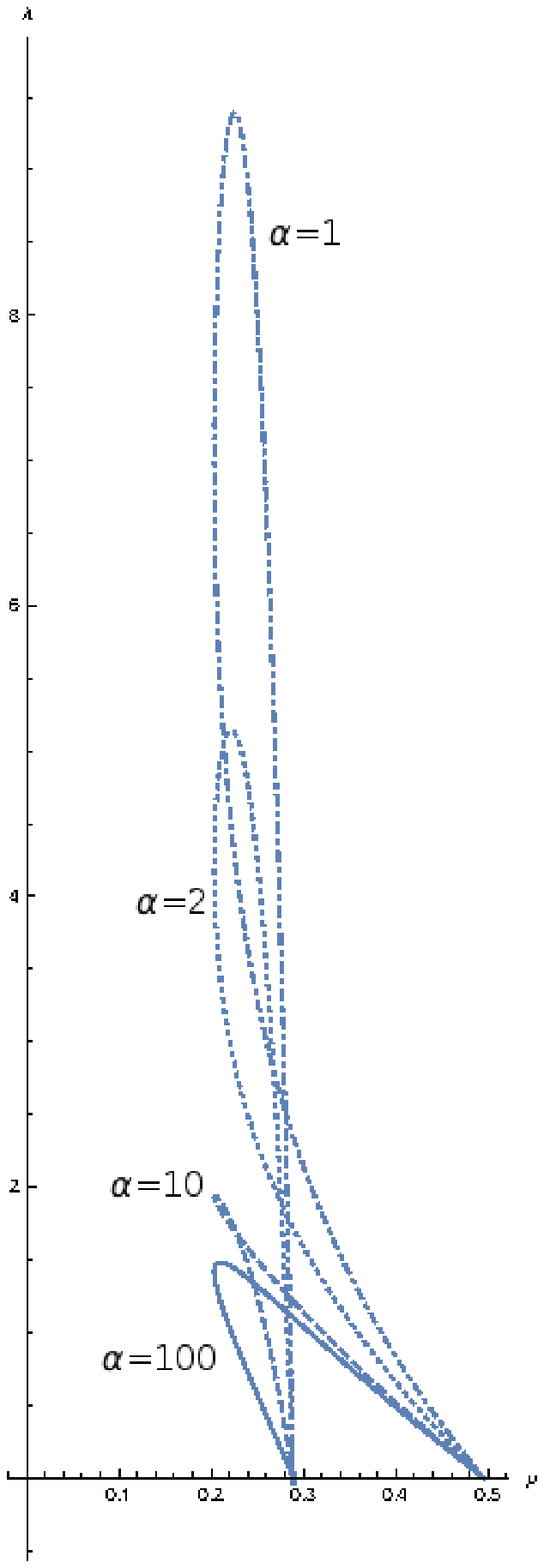}
 \caption{The Hopf curves for the system \eqref{eq:Atri3Da}--\eqref{eq:Atri3Dc} when  $T(c)=\frac{\alpha c^2}{1+\alpha c^2}$ for $\alpha=1, 2, 10, 100$.}
\label{fig:HopfCurvesAlphaVariesNewHillFn}
\end{center}
\end{figure}

In Figure \ref{fig:HopfCurvesAlphaVariesNewHillFn}  we plot the Hopf curves of the system \eqref{eq:Atri3Da}--\eqref{eq:Atri3Dc} for the Hill function $T=\frac{\alpha c^2}{1+\alpha c^2}$, selecting again different values of $\alpha$.
Comparing  Figure \ref{fig:HopfCurvesAlphaVariesNewHillFn} with Figure \ref{fig:HopfCurvesAlphaVaries} for $T(c)=\frac{\alpha c}{1+\alpha c}$ we see that the system has the same qualitative behaviour as $\alpha$ decreases; for a certain value of $\alpha$ the system develops a cusp and for values of $\alpha$ smaller than the value at the cusp the system has a "bow-tie".

The parametric method we developed above is very useful as we can easily plot the Hopf curve, and the other important curves of the bifurcation diagram, not only for any Hill function, but also for any other functional form of $T$. In this work, we have considered functions of $T$ that satisfy $T(0)=0$ and saturate to a constant value of $T$ as $c$ gets large, but the same method could be used for a $T(c)$ with different qualitative characteristics, as could possibly arise from experimental work.
 
\subsection{Bifurcation diagrams and limit cycles}
\begin{figure}[ht] 
  \begin{subfigure}[b]{0.5\linewidth}
    \centering
    \includegraphics[width=1\linewidth]{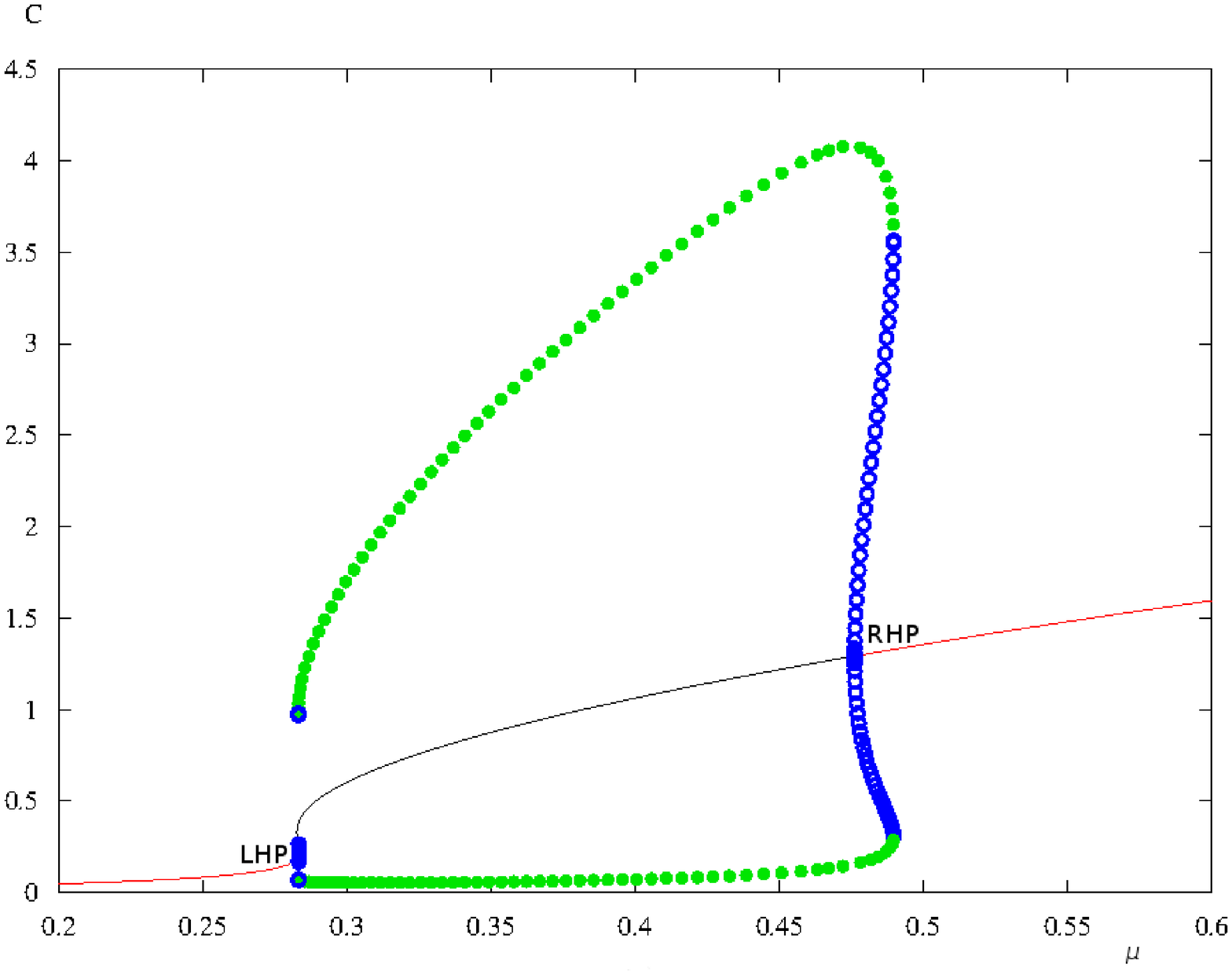}
    \caption{$\lambda=0.1$} 
  \label{fig:AtriMechsLam0Pt1}
    \vspace{4ex}
  \end{subfigure}
  \begin{subfigure}[b]{0.5\linewidth}
    \centering
    \includegraphics[width=1\linewidth]{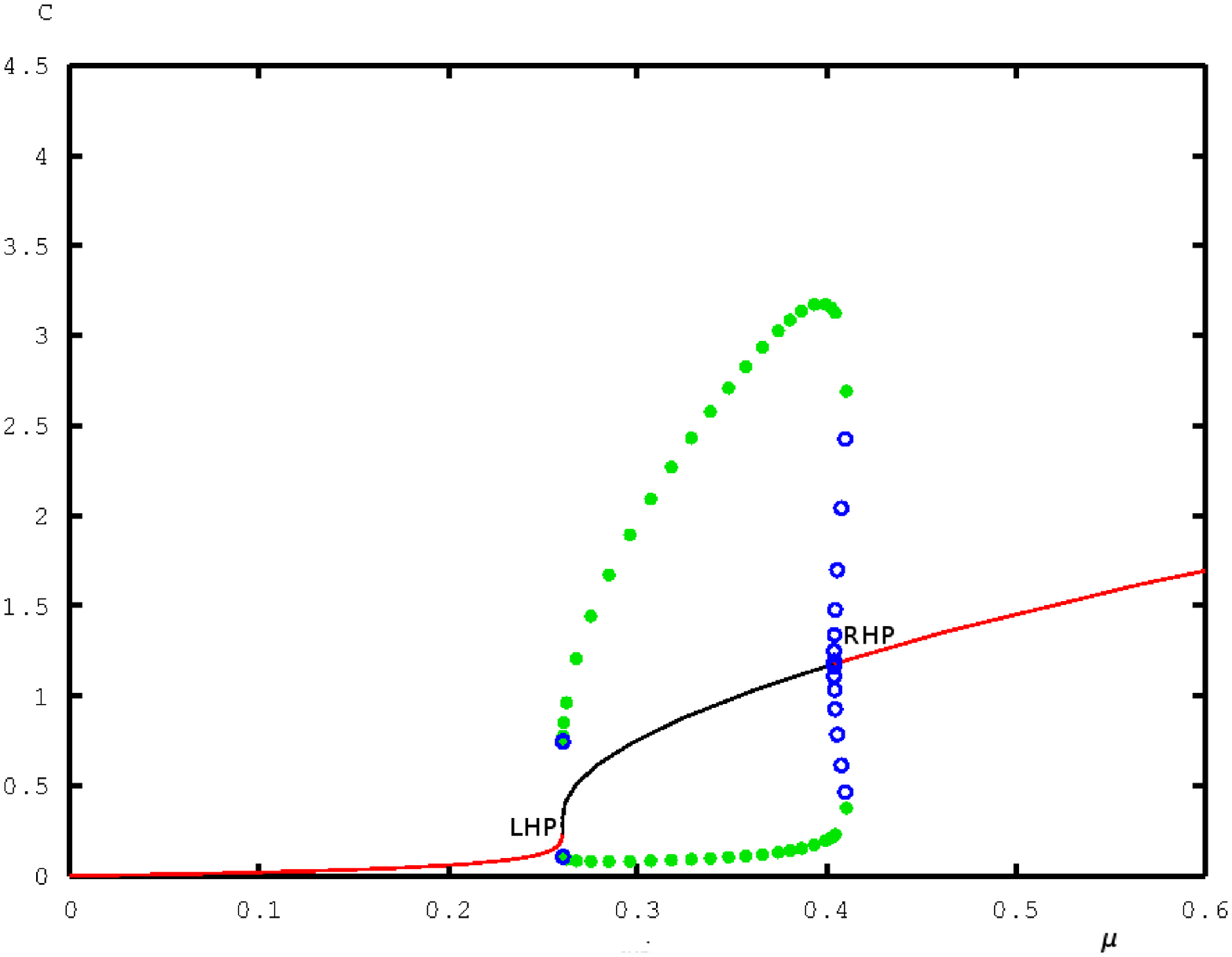}
    \caption{$\lambda=0.5$} 
   \label{fig:AtriMechsLam0Pt5}
    \vspace{4ex}
  \end{subfigure} 
  \begin{subfigure}[b]{0.5\linewidth}
    \centering
    \includegraphics[width=1\linewidth]{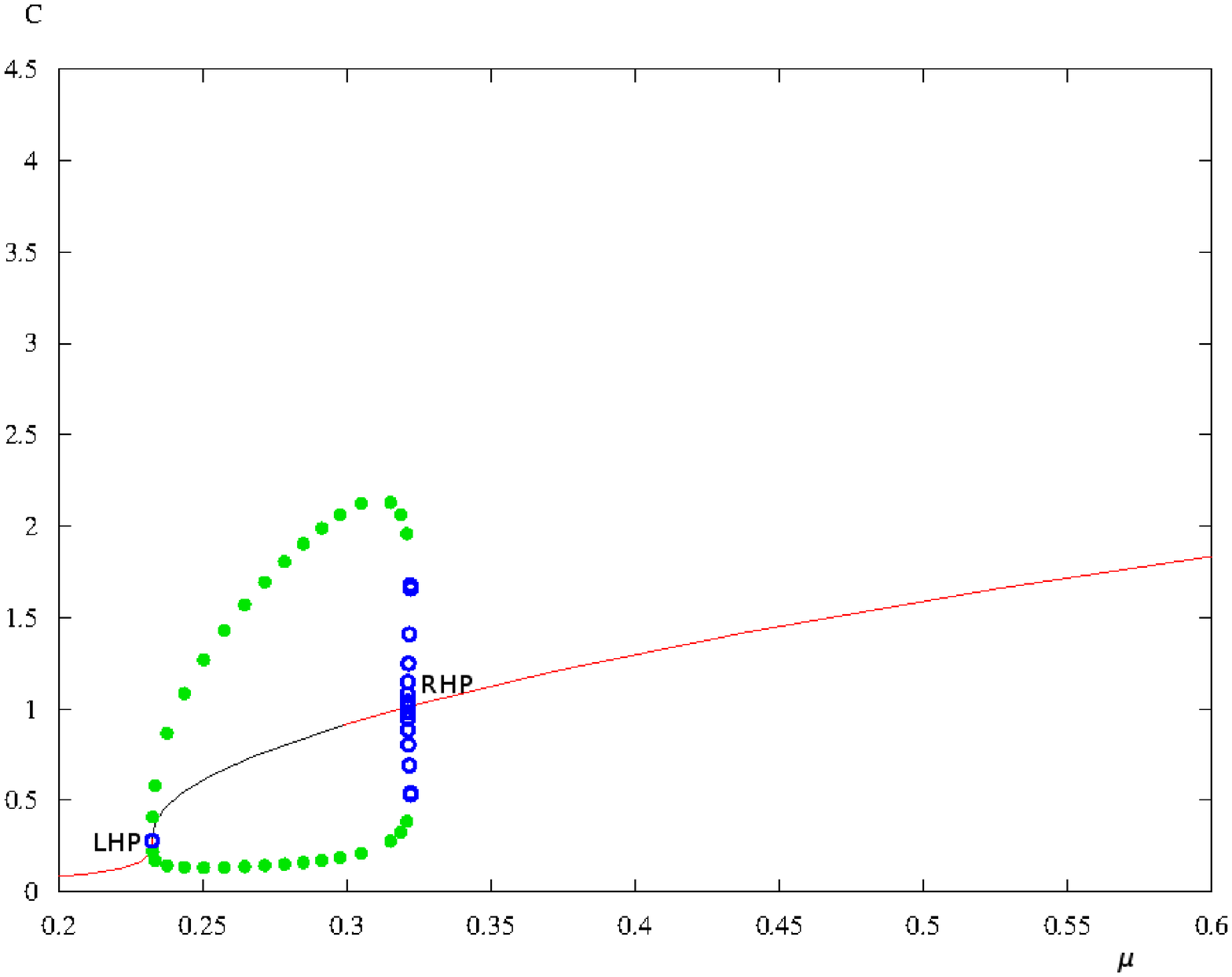}
    \caption{$\lambda=1$} 
  \label{fig:AtriMechsLam1} 
  \end{subfigure}
  \begin{subfigure}[b]{0.5\linewidth}
    \centering
     \includegraphics[width=1\linewidth]{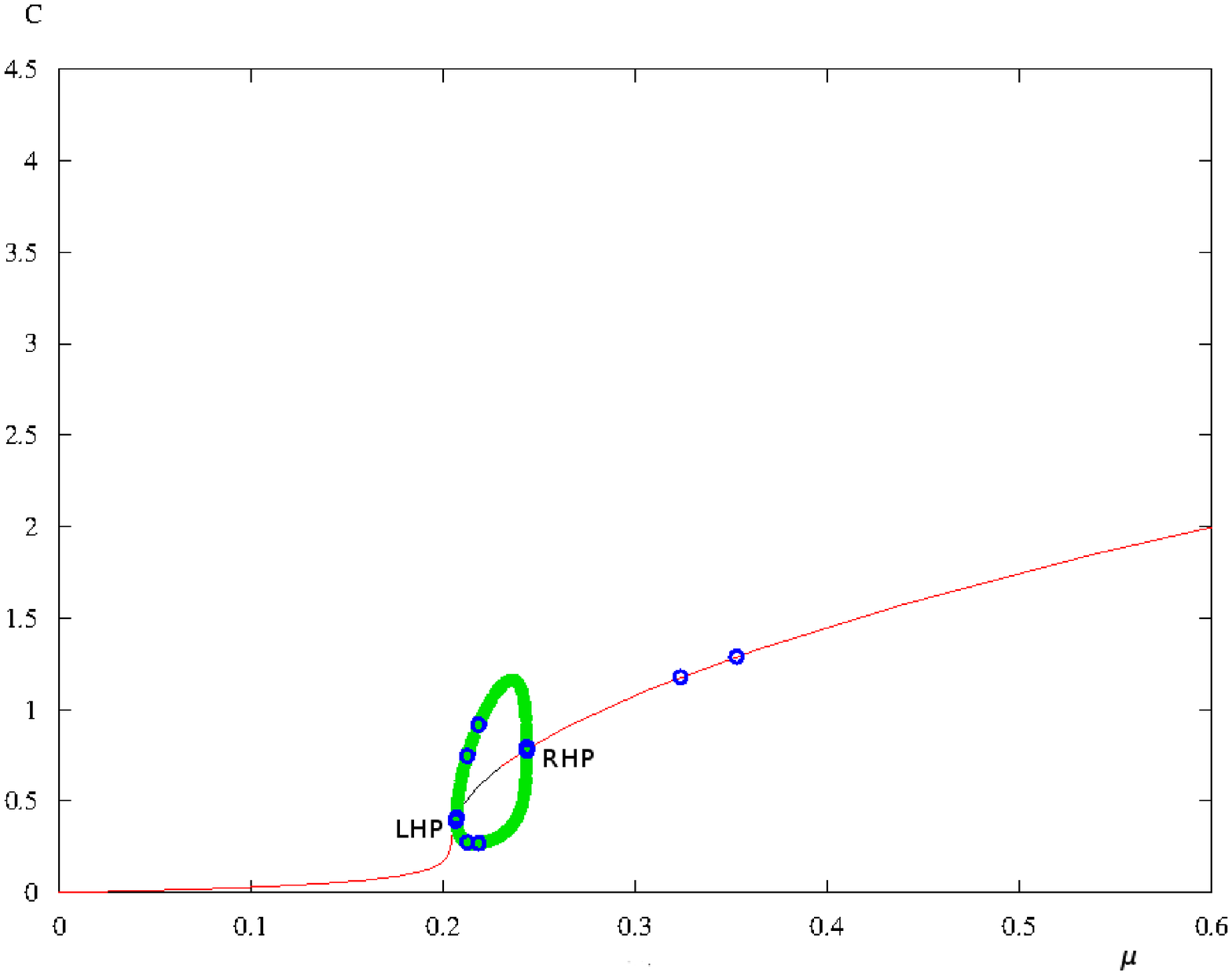}
    \caption{$\lambda=1.5$} 
   \label{fig:AtriMechsLam1Pt5}
  \end{subfigure} 
  \caption{Bifurcation diagrams for the system \eqref{eq:Atri3Da}--\eqref{eq:Atri3Dc} when $T(c)=\frac{10 c}{1+10 c}$, as $\mu$ is increased and for selected values of $\lambda$. The LHP and the RHP are indicated. The stable limit cycles are represented by solid circles and the unstable limit cycles by open circles (respectively with green and blue colour on the online version of the article): (a) $\lambda=0.1$ (b) $\lambda=0.5$ (c) $\lambda=1$ (d) $\lambda=1.5$ (note that there are two extraneous open circles in the stable branch of the equilibrium curve, a numerical artifact of XPPAUT).}
\label{fig:AtriMechs1}
\end{figure}

\begin{figure}[ht] 
  \begin{subfigure}[b]{0.5\linewidth}
    \centering
    \includegraphics[width=0.99\linewidth]{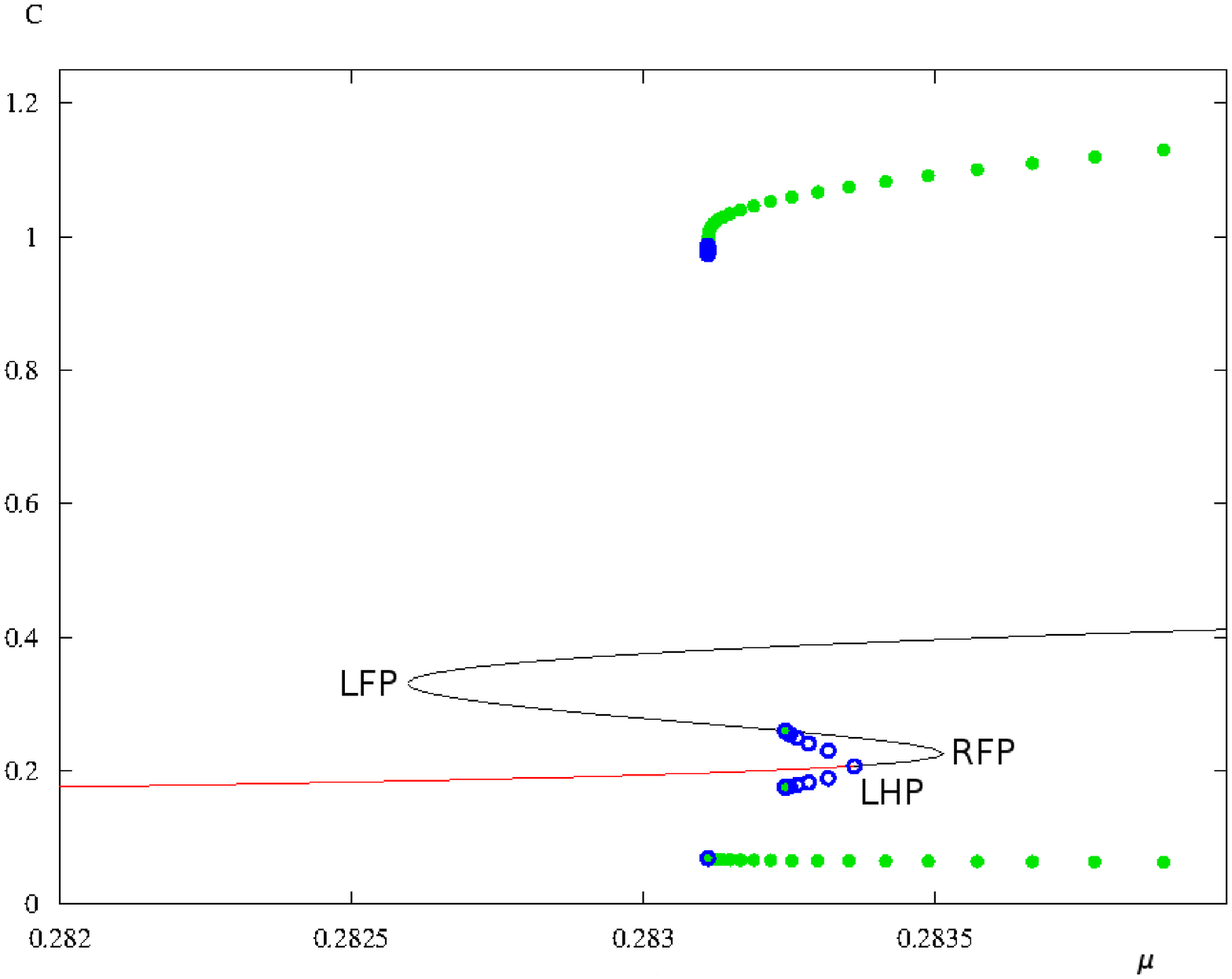}
    \caption{$\lambda=0.1$} 
  \label{fig:AtriMechsLam0Pt1Zoom}
    \vspace{4ex}
  \end{subfigure}
  \begin{subfigure}[b]{0.5\linewidth}
    \centering
    \includegraphics[width=0.99\linewidth]{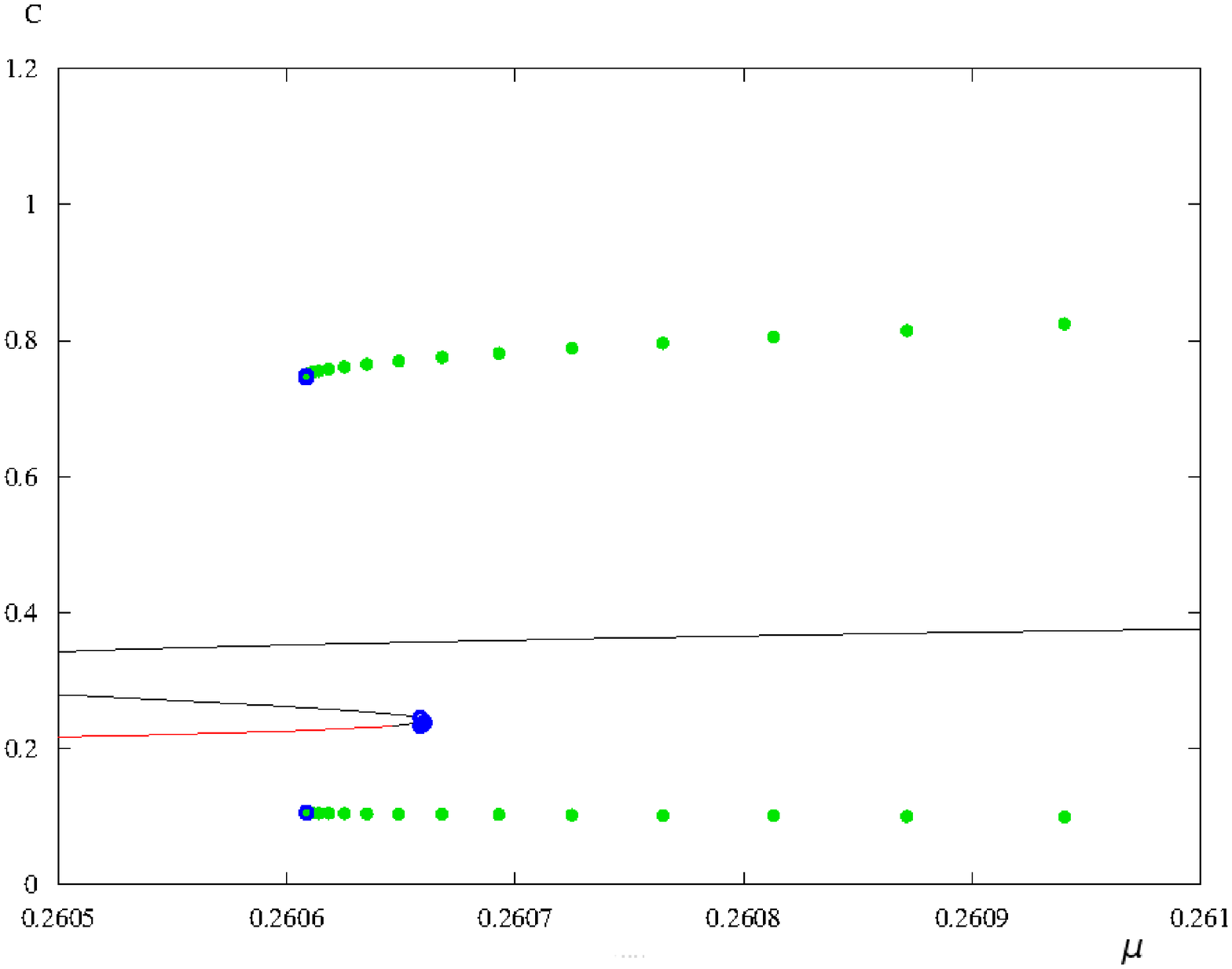}
    \caption{$\lambda=0.5$} 
   \label{fig:AtriMechsLam0Pt5Zoom}
    \vspace{4ex}
  \end{subfigure} 
  \begin{subfigure}[b]{1\linewidth}
    \centering
    \includegraphics[width=0.5\linewidth]{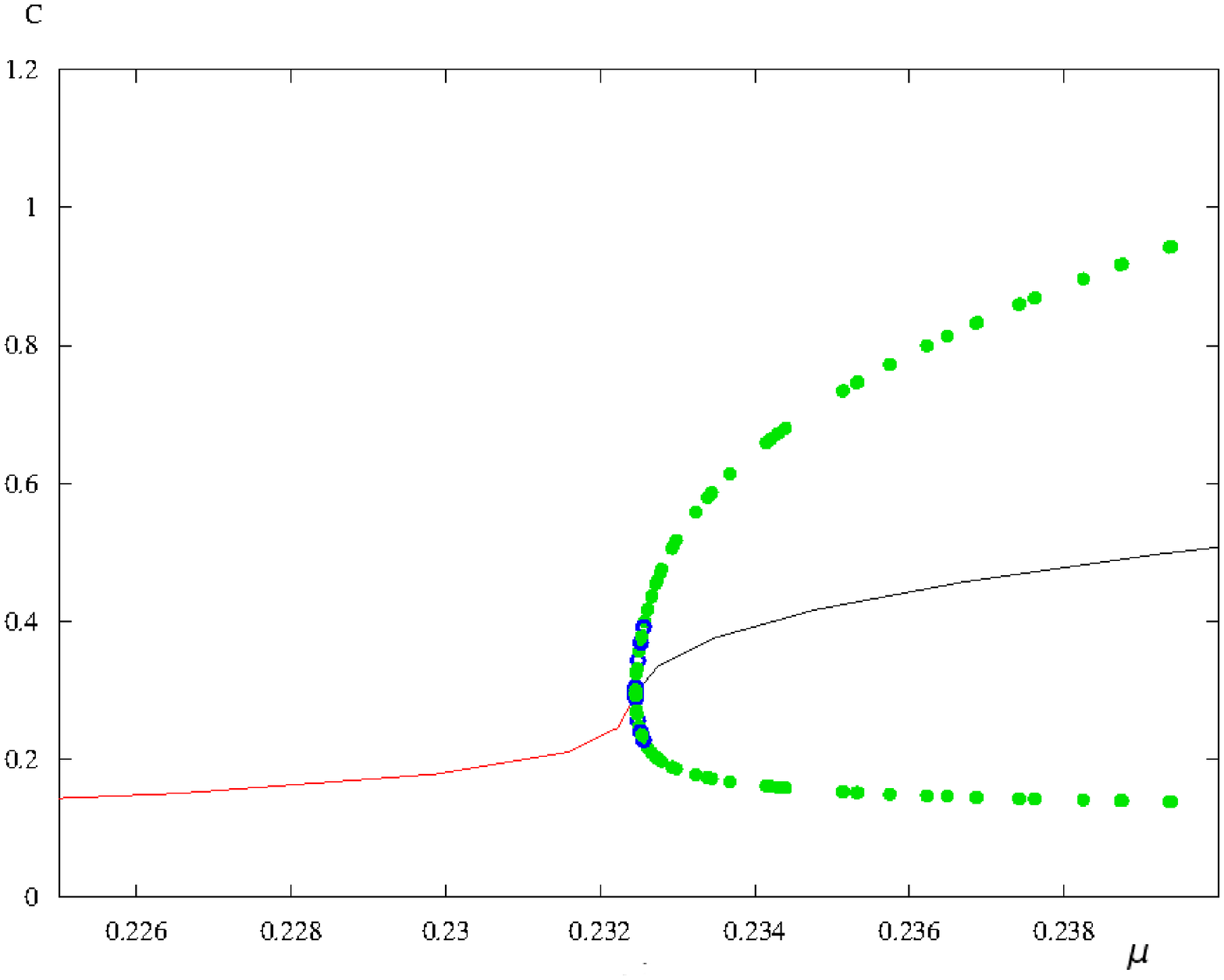}
    \caption{$\lambda=1$} 
  \label{fig:AtriMechsLam1Zoom} 
  \end{subfigure}
  \caption{Zooming in on the left Hopf point of the system \eqref{eq:Atri3Da}--\eqref{eq:Atri3Dc} when $T(c)=\frac{10 c}{1+10 c}$ for: (a) $\lambda=0.1$ (b) $\lambda=0.5$ (c) $\lambda=1$, as $\mu$ increases.}
\label{fig:AtriMechs1Zoom}
\end{figure}

We now use AUTO (as implemented in XPPAUT) to determine the amplitude and frequency of oscillations (limit cycles) of the system \eqref{eq:Atri3Da}--\eqref{eq:Atri3Dc} when $T(c)=10c/(1+10c)$. We plot one-parameter bifurcation diagrams as $\mu$ is increased, for different values of $\lambda$, in Figure \ref{fig:AtriMechs1}. We see that as $\lambda$ increases the amplitude of oscillations decreases until the oscillations vanish (close to) $\lambda=\lambda_{\rm  max}=1.68632$, which agrees with the linear analysis (see Figure \ref{fig:HopfCurveFold}). We also observe that for $\lambda=0.1, 0.5$ and $1$, respectively in Figures \ref{fig:AtriMechsLam0Pt1}, \ref{fig:AtriMechsLam0Pt5}, and \ref{fig:AtriMechsLam1},  there are both stable and unstable limit cycles (solid and open circles, respectively), and the right Hopf point is subcritical. The $\mu$-range of unstable limit cycles decreases as $\lambda$ increases. For $\lambda=1.5$ (Figure \ref{fig:AtriMechsLam1Pt5}) there are only stable limit cycles and the right Hopf point has become supercritical. In Figure \ref{fig:AtriMechs1Zoom} we zoom in on the left Hopf point for $\lambda=0.1, 0.5$ and $1$. We see that the branches of unstable limit cycles coming out of the left Hopf point shorten as 
$\lambda$ increases. For $\lambda=1$ the limit cycle branches coming out of the left and right Hopf point join up, with only stable limit cycles present. 

\begin{figure}[ht] 
  \begin{subfigure}[b]{0.5\linewidth}
    \centering
    \includegraphics[width=0.99\linewidth]{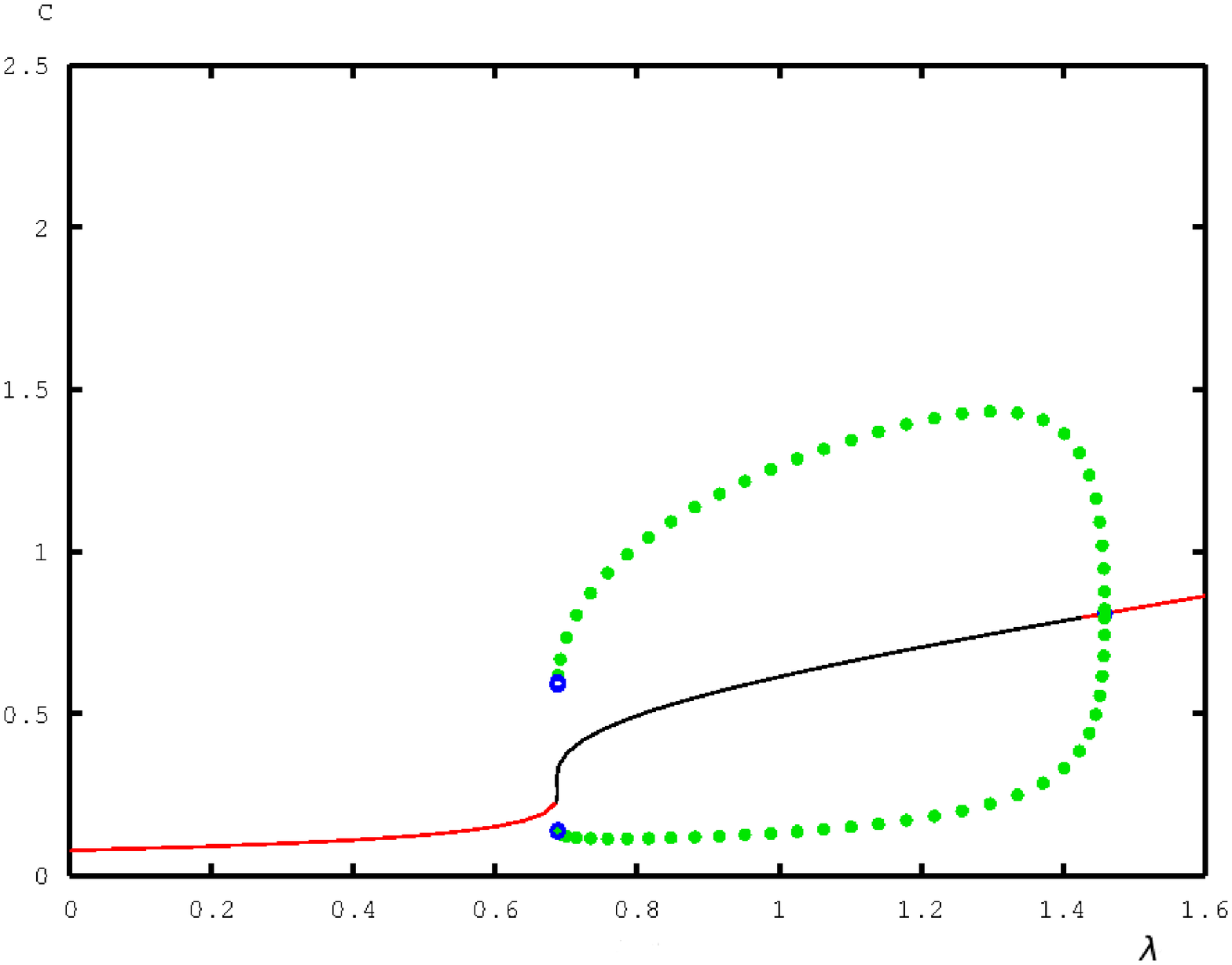}
    \caption{$\mu=0.25$} 
  \label{fig:AtriMechsMu0Pt25}
    \vspace{4ex}
  \end{subfigure}
  \begin{subfigure}[b]{0.5\linewidth}
    \centering
    \includegraphics[width=0.99\linewidth]{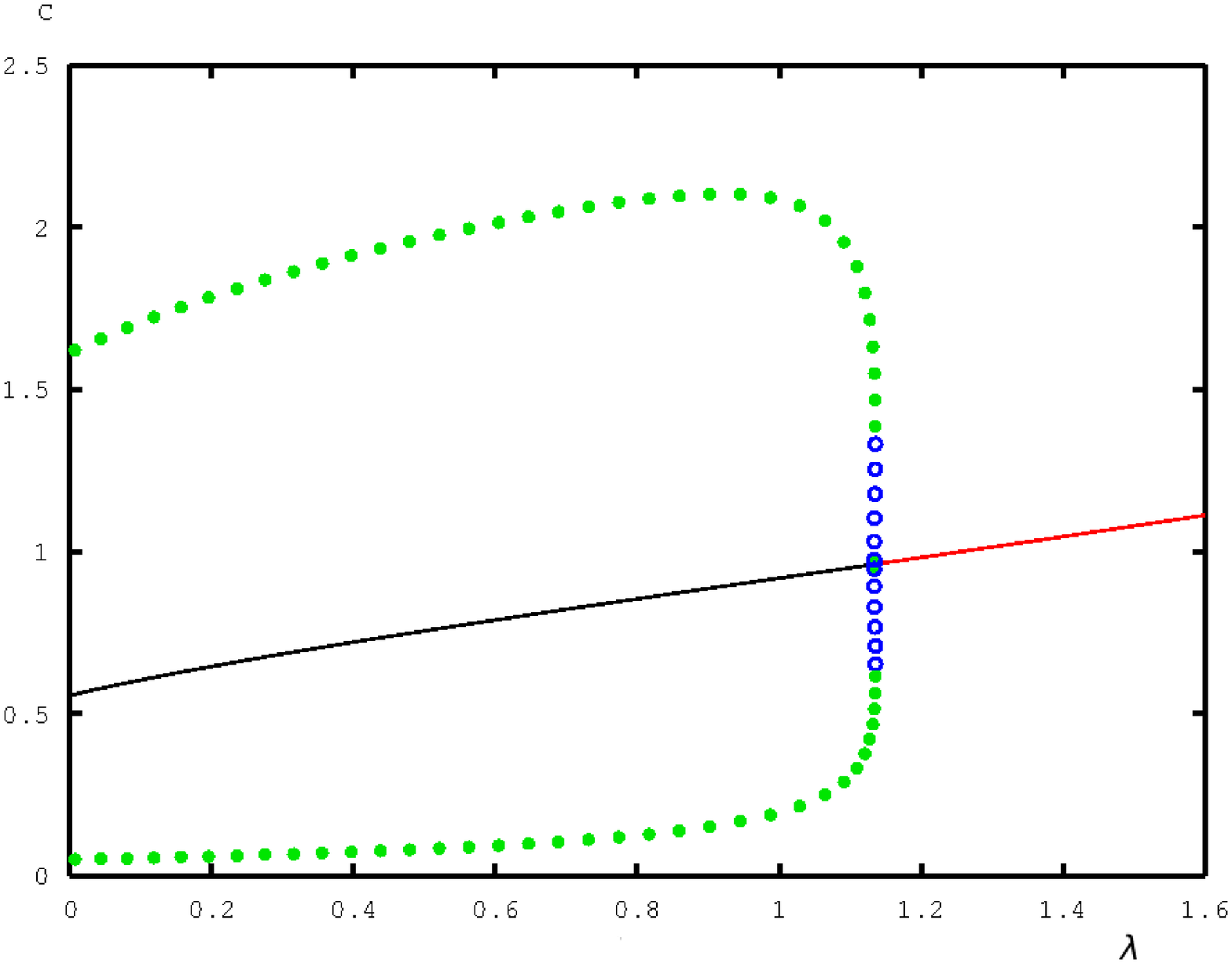}
    \caption{$\mu=0.3$} 
  \label{fig:AtriMechsMu0Pt3} 
    \vspace{4ex}
  \end{subfigure} 
  \caption{Bifurcation diagrams for the system \eqref{eq:Atri3Da}--\eqref{eq:Atri3Dc} when  $T(c)=\frac{10 c}{1+10 c}$, showing the amplitude of oscillations as $\lambda$ is increased, for: (a) $\mu=0.25$ (b) $\mu=0.3$. The stable limit cycles are represented by solid circles and the unstable limit cycles by open circles-respectively with green and blue colour on the online version of the article.}
\label{fig:AtriMechsMu}
\end{figure}

In Figure \ref{fig:AtriMechsMu} we plot one-parameter bifurcation diagrams as $\lambda$ is increased, for two selected values of $\mu$. For $\mu=0.25$the Atri system ($\lambda=0$) has no oscillations but for the mechanochemical model, as $\lambda$ is increased, we see in Figure \ref{fig:AtriMechsMu0Pt25} that stable limit cycles (this, ties up with Figure \ref{fig:HopfCurveFold}). For $\mu=0.3$ the Atri system has a stable limit cycle. As $\lambda$ is increased, we see in Figure \ref{fig:AtriMechsMu0Pt3} that stable and unstable limit cycles emerge for a finite range of $\lambda$ values.

\begin{figure}[ht] 
  \begin{subfigure}[b]{0.5\linewidth}
    \centering
    \includegraphics[width=0.99\linewidth]{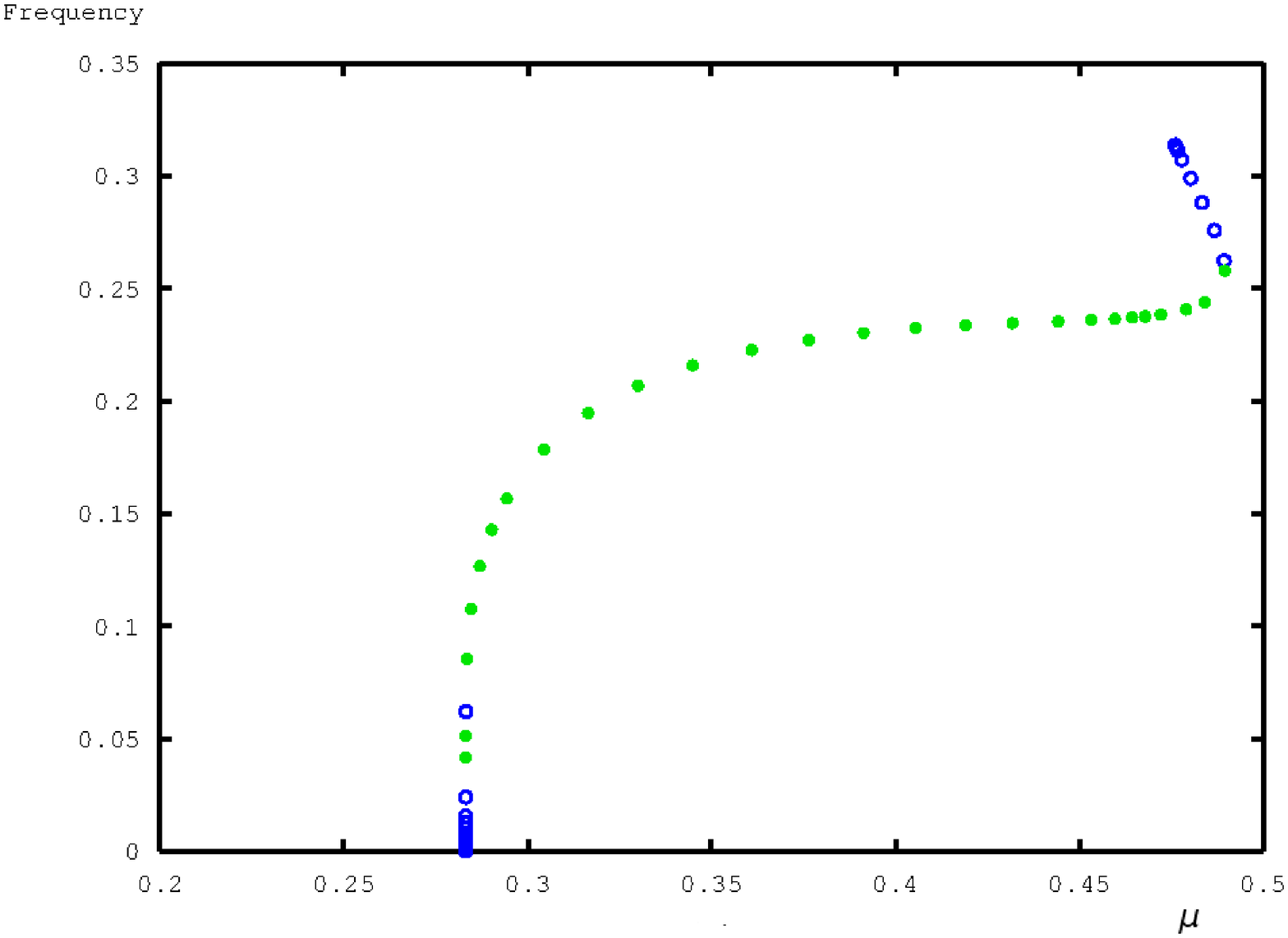}
    \caption{$\lambda=0.1$} 
  \label{fig:Lam0Pt1Freq}
    \vspace{4ex}
  \end{subfigure}
  \begin{subfigure}[b]{0.5\linewidth}
    \centering
    \includegraphics[width=0.99\linewidth]{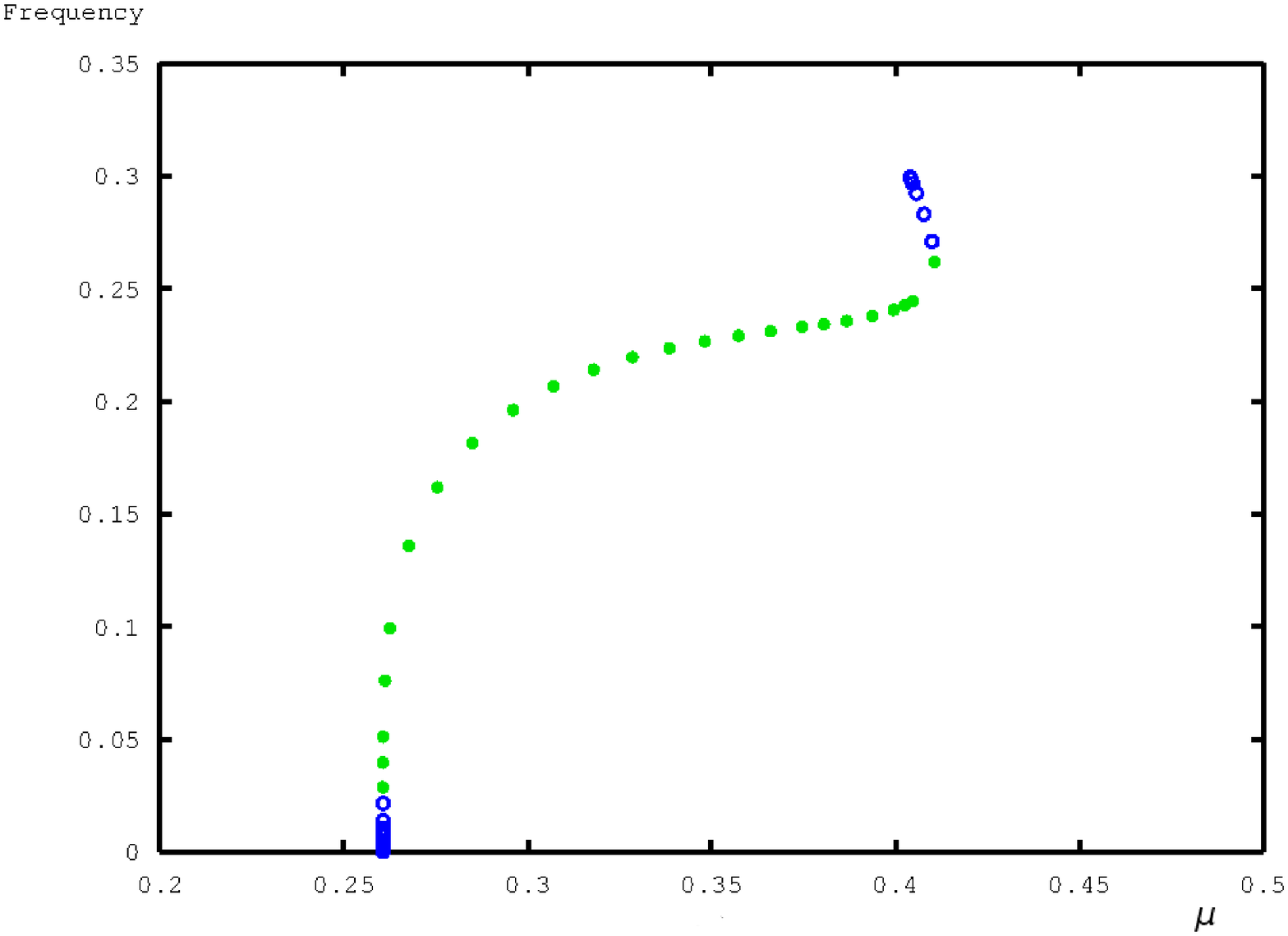}
    \caption{$\lambda=0.5$} 
   \label{fig:Lam0Pt5Freq}
    \vspace{4ex}
  \end{subfigure} 
  \begin{subfigure}[b]{0.5\linewidth}
    \centering
    \includegraphics[width=0.99\linewidth]{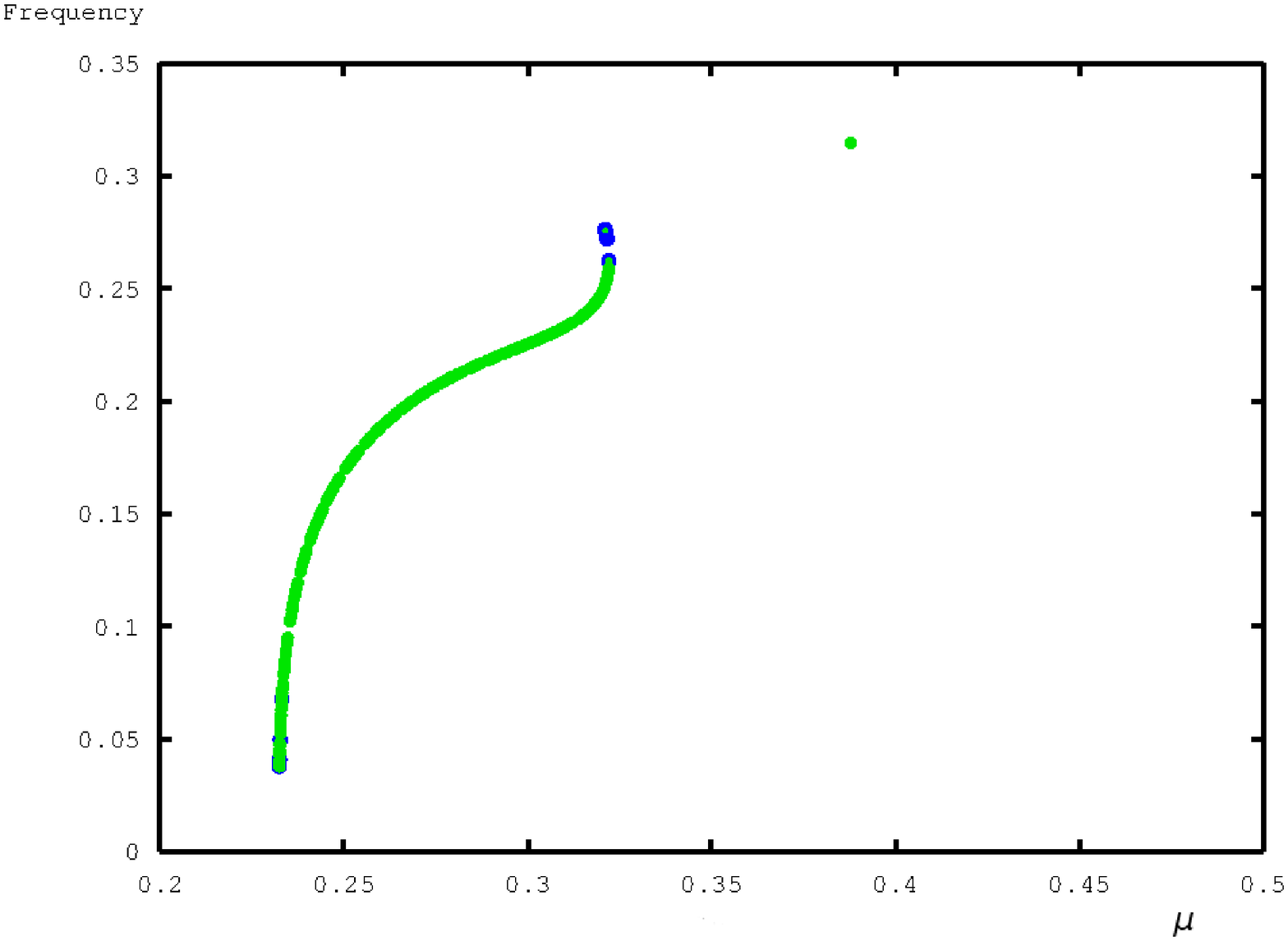}
    \caption{$\lambda=1$} 
  \label{fig:Lam1Freq} 
  \end{subfigure}
  \begin{subfigure}[b]{0.5\linewidth}
    \centering
    \includegraphics[width=0.99\linewidth]{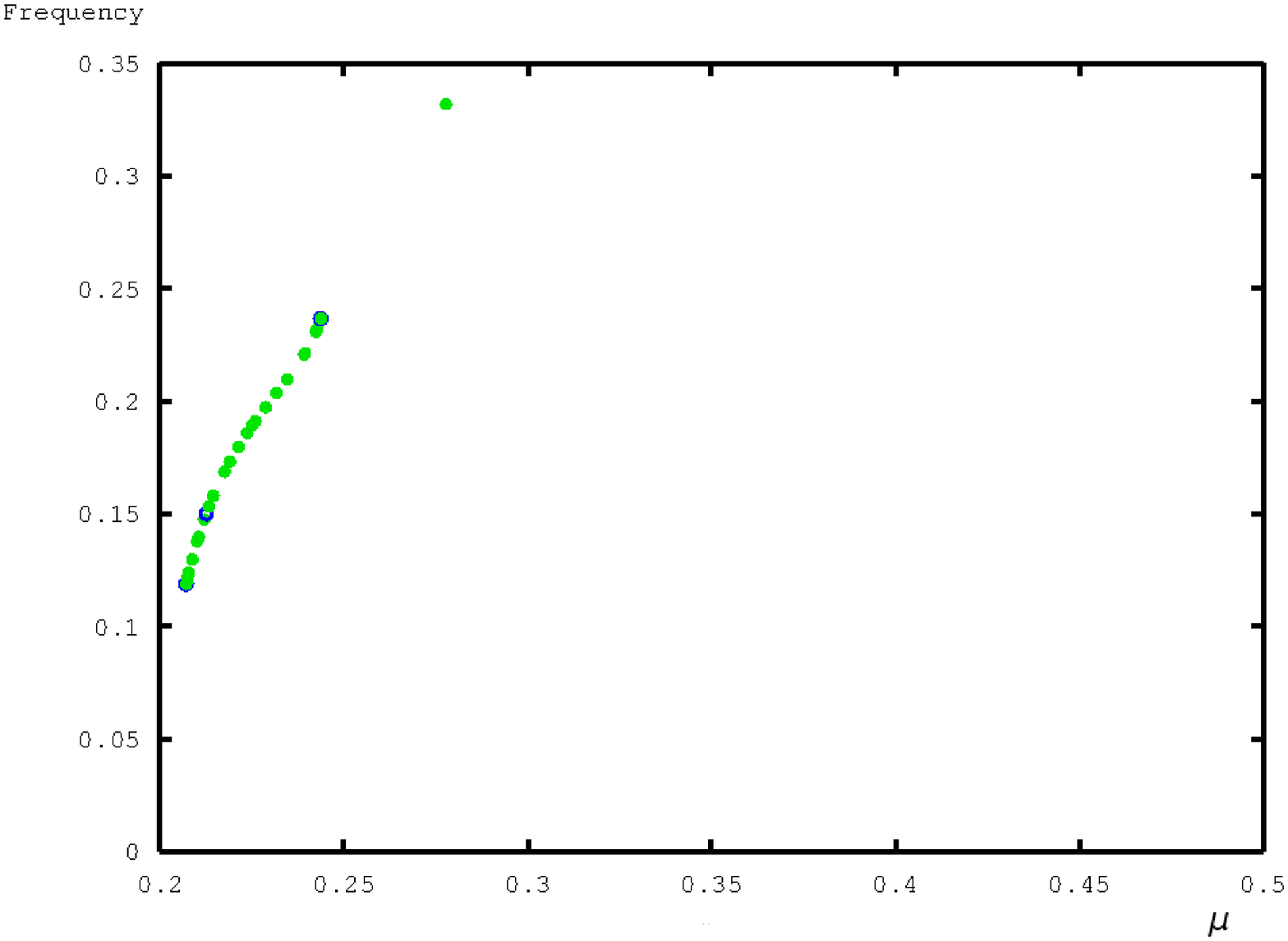}
    \caption{$\lambda=1.5$} 
   \label{fig:Lam1Pt5Freq}
  \end{subfigure} 
  \caption{Frequency of oscillations for the system \eqref{eq:Atri3Da}--\eqref{eq:Atri3Dc} when  $T(c)=\frac{10 c}{1+10 c}$ as a function of $\mu$ and selected values of $\lambda$: (a) $\lambda=0.1$ (b) $\lambda=0.5$ (c) $\lambda=1$ (d) $\lambda=1.5$. The stable limit cycles are represented by solid circles and the unstable limit cycles by open circles-respectively with green and blue colour on the online version of the article (note that there is one extraneous solid (green) circle in Figure  \ref{fig:Lam1Freq}, and in Figure \ref{fig:Lam1Pt5Freq}, both numerical artifacts of the XPPAUT software. (This is also confirmed by looking at the range of oscillations in Figures \ref{fig:AtriMechsLam1} and \ref{fig:AtriMechsLam1Pt5}).}
\label{fig:AtriMechsFreq}
\end{figure}

In Figure \ref{fig:AtriMechsFreq} we plot the frequency of the limit cycles for $\lambda=0.1, 0.5, 1, 1.5$. For $\lambda=0.1, 0.5, 1$ the frequency rapidly increases close to the LHP and the RHP and there is an `intermediate' region where the frequency varies slowly with $mu$. This is qualitatively similar to the $\lambda=0$ case (Atri model)-cf. Figure \ref{fig:AtriBifnFrequency}. The intermediate region gets smaller as $\lambda$ increases, and for $\lambda=1.5$ (Figure \ref{fig:Lam1Pt5Freq}) the latter region seems to have vanished. These qualitative changes in the oscillation frequency as $\lambda$ is increased seem to correlate with the qualitative changes in the oscillation amplitude, as seen in Figure \eqref{fig:AtriMechs1}.
Summarising, for any value of $\mu$ and $\lambda$ we can determine the approximate \emph{range of oscillations} using the parametric expressions \eqref{eq:MUvsC} and \eqref{eq:LAMBDAvsC}, and we can use XPPAUT to also obtain the \emph{amplitude and the frequency of oscillations}. 

\begin{figure}[ht] 
  \begin{subfigure}[b]{0.5\linewidth}
    \centering
    \includegraphics[width=0.95\linewidth]{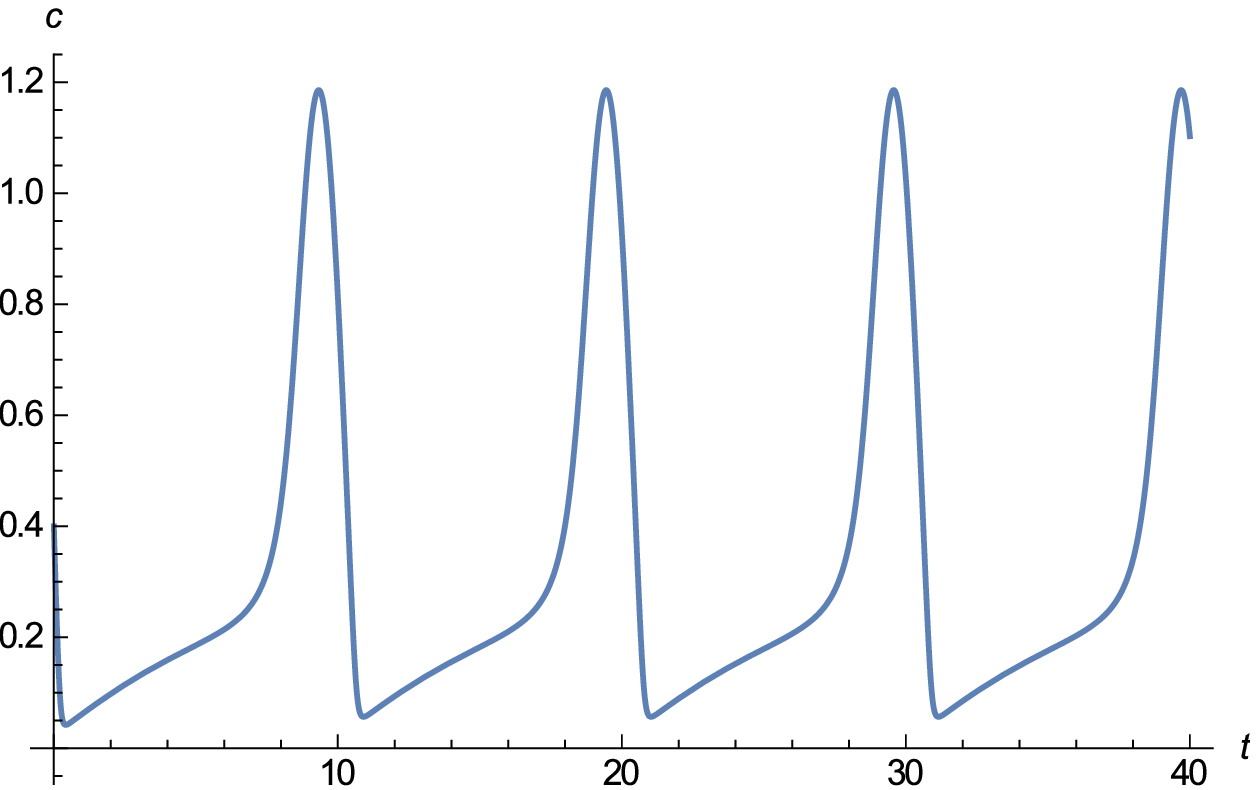}
    \caption{$\lambda=0$} 
  \label{fig:1a}		
    \vspace{4ex}
  \end{subfigure}
  \begin{subfigure}[b]{0.5\linewidth}
    \centering
    \includegraphics[width=0.95\linewidth]{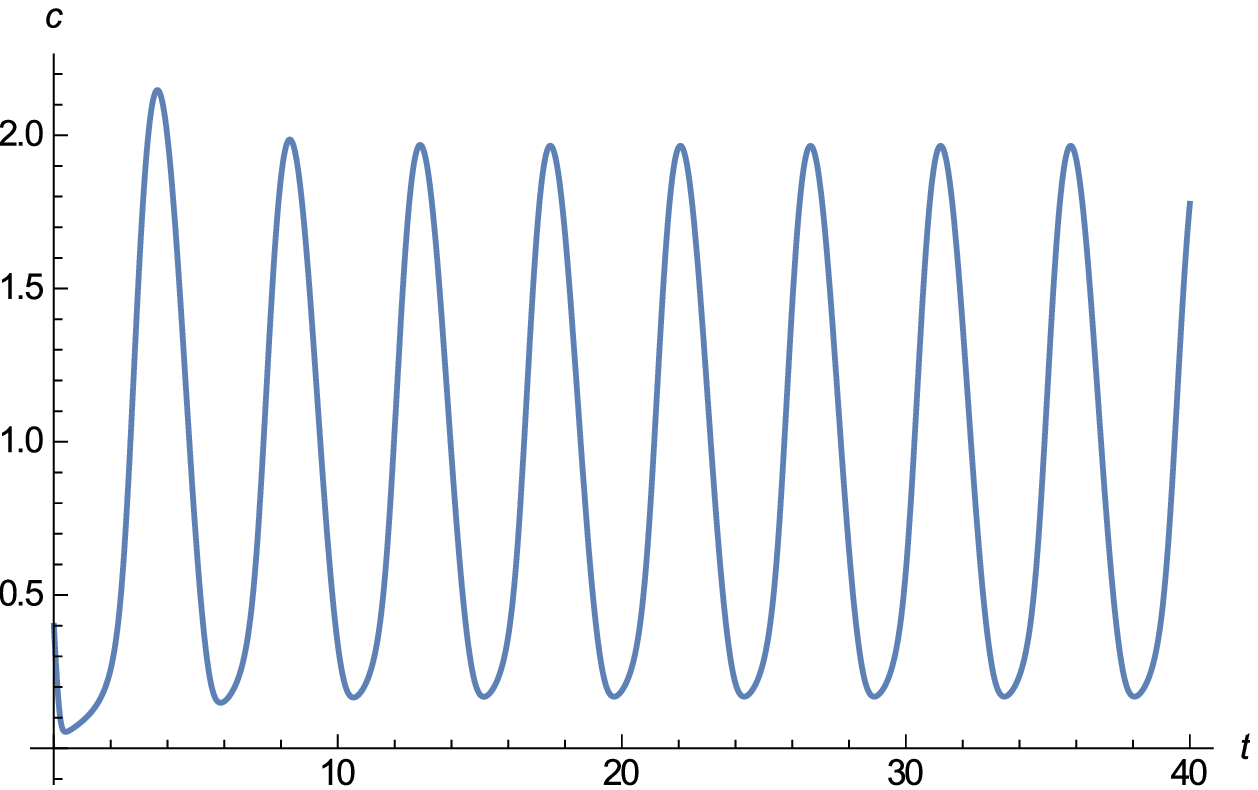}
		\caption{$\lambda=1$}\label{fig:1b}	
    \vspace{4ex}
  \end{subfigure} 
  \begin{subfigure}[b]{1\linewidth}
    \centering
    \includegraphics[width=0.5\linewidth]{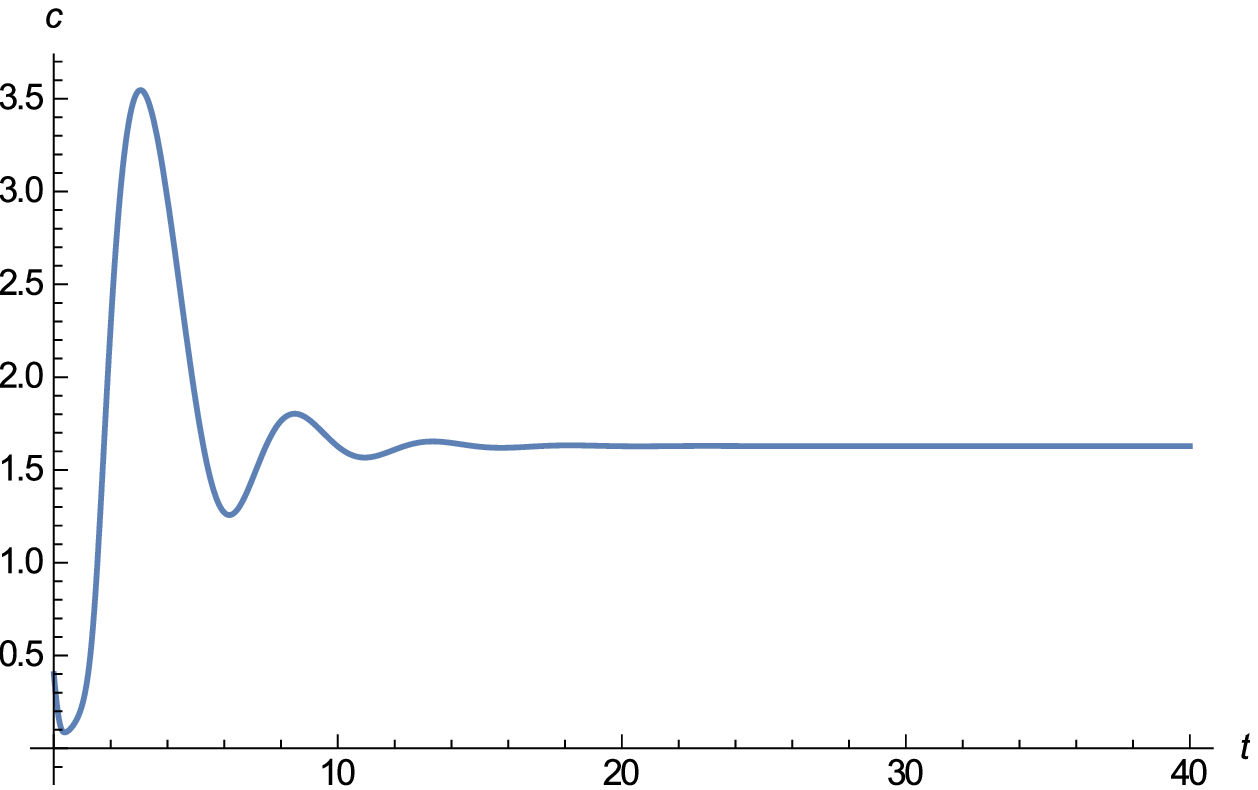}
		\caption{$\lambda=3$}\label{fig:1c}
  \end{subfigure}
  \caption{Evolution of $c(t)$ with time, solving the system \eqref{eq:Atri3Da}--\eqref{eq:Atri3Dc} when $T(c)=\frac{10 c}{1+10 c}$, $\mu=0.2894$: (a) $\lambda=0$ (Atri model) (b) $\lambda=1$ (the frequency and amplitude of oscillations has increased)  (c) $\lambda=3$ (oscillations are suppressed).}
\label{fig:1}
\end{figure}

In Figure \ref{fig:1} we plot the evolution of $c(t)$, solving \eqref{eq:Atri3Da}--\eqref{eq:Atri3Dc} for $\mu=0.2894$ and selected values of $\lambda$; as expected from the bifurcation diagrams above, the oscillations  are suppressed when $\lambda$ is sufficiently increased.


\subsection{Asymptotic analysis}
We now develop an asymptotic analysis for the mechanochemical model \eqref{eq:Atri3Da}--\eqref{eq:Atri3Dc}, similar to that for the Atri model, assuming that $c$ is a fast variable compared to $h$ and $\theta$.
We divide equation \eqref{eq:Atri3Da} by $K_1$ and on the left hand side we replace $1/K_1$ by $\epsilon$, a small parameter. Therefore, the ODE system becomes
\begin{align}
\label{eq:Atri3Da-Asym}
\epsilon \frac{dc}{dt}&=\mu h\frac{b+c}{1+c}-\frac{\Gamma}{K_1}\frac{c}{K+c}+\frac{\lambda}{K_1} \theta=R_1/K_1=F_1,\\
\label{eq:Atri3Db-Asym}
\frac{d\theta}{dt}&=-\theta+T(c)=R_2,\\
\label{eq:Atri3Dc-Asym}
\frac{dh}{dt}&=\frac{1}{1+c^2}-h=R_3.
\end{align}
For the asymptotic analysis we let $\epsilon \to 0$ and we also assume that $\frac{\Gamma}{K_1}$ and $\Lambda= \frac{\lambda}{K_1}$ are O(1).
\subsubsection{Fast system}
The transformation $t=\epsilon \tau$ converts the system \eqref{eq:Atri3Da-Asym}--\eqref{eq:Atri3Dc-Asym} into the fast system
\begin{align}
\label{eq:Atri3Da-Asym-Fast}
\frac{dc}{d\tau}&=\mu h\frac{b+c}{1+c}-\frac{\Gamma}{K_1}\frac{c}{K+c}+\Lambda\theta=F_1,\\
\label{eq:Atri3Db-Asym-Fast}
\frac{d\theta}{d\tau}&=\epsilon(-\theta+T(c))=\epsilon R_2,\\
\label{eq:Atri3Dc-Asym-Fast}
\frac{dh}{d\tau}&=\epsilon\left(\frac{1}{1+c^2}-h\right)=\epsilon R_3.
\end{align}
Inserting $c=c_0(\tau)+O(\epsilon)$, $\theta=\theta_0+O(\epsilon)$, $h=h_0(\tau)+O(\epsilon)$ in \eqref{eq:Atri3Da-Asym-Fast}--\eqref{eq:Atri3Dc-Asym-Fast} we obtain, at leading order,
\begin{align}
\label{eq:Atri3Da-Asym-Fast-LO}
\frac{dc_0}{d\tau}&=\mu h_0\frac{b+c_0}{1+c_0}-\frac{\Gamma}{K_1}\frac{c_0}{K+c_0}+\Lambda\theta_0\\
\label{eq:Atri3Db-Asym-Fast-LO}
\frac{d\theta_0}{d\tau}&=0\implies \theta_0=\theta_{00}={\rm constant,}\,\,\,\frac{dh_0}{d\tau}=0\implies h_0=h_{00}={\rm constant}.
\end{align}
Therefore, from \eqref{eq:Atri3Db-Asym-Fast-LO} the slow variables $h$ and $\theta$ are constant at leading order and only the fast variable $c$ varies according to equation \eqref{eq:Atri3Da-Asym-Fast-LO}. The planes $h_0=h_{00}={\rm constant}$ and  $\theta_0=\theta_{00}={\rm constant}$ fill the whole of $R^3$ and constitute the \emph{fast foliation} of the space, describing the evolution of the system on the fast timescale. Since the intersection of any two planes $h$=constant and $\theta$=constant is a line the fast motion of the system is along straight lines.

\subsubsection{Slow system and slow manifold}
\begin{figure}	
	\centering
	\begin{subfigure}[t]{2.3in}
		\centering
		\includegraphics[width=2.3in]{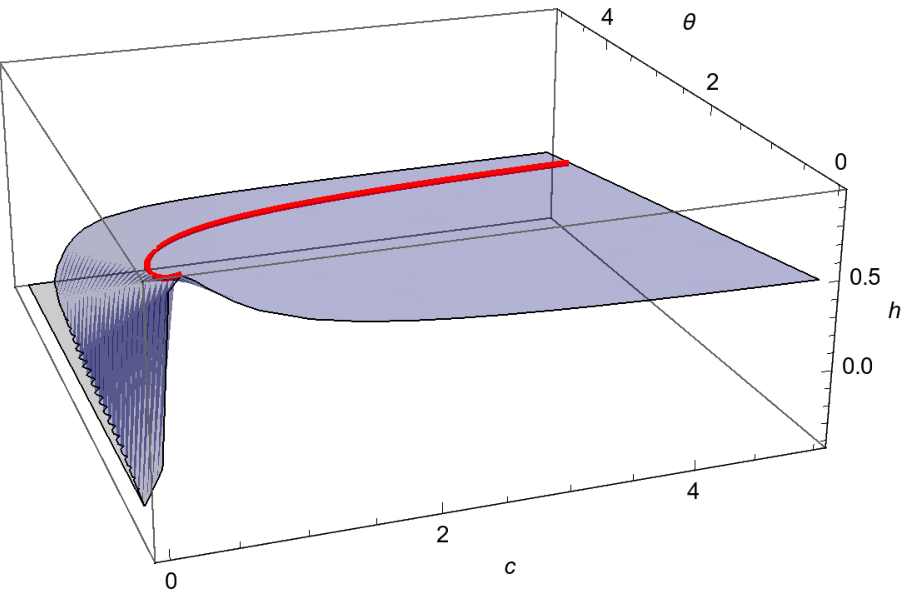}
		\caption{}\label{fig:SMandFOLD-BetterView}		
	\end{subfigure}
	\quad
	\begin{subfigure}[t]{2in}
		\centering
		\includegraphics[width=1.8in]{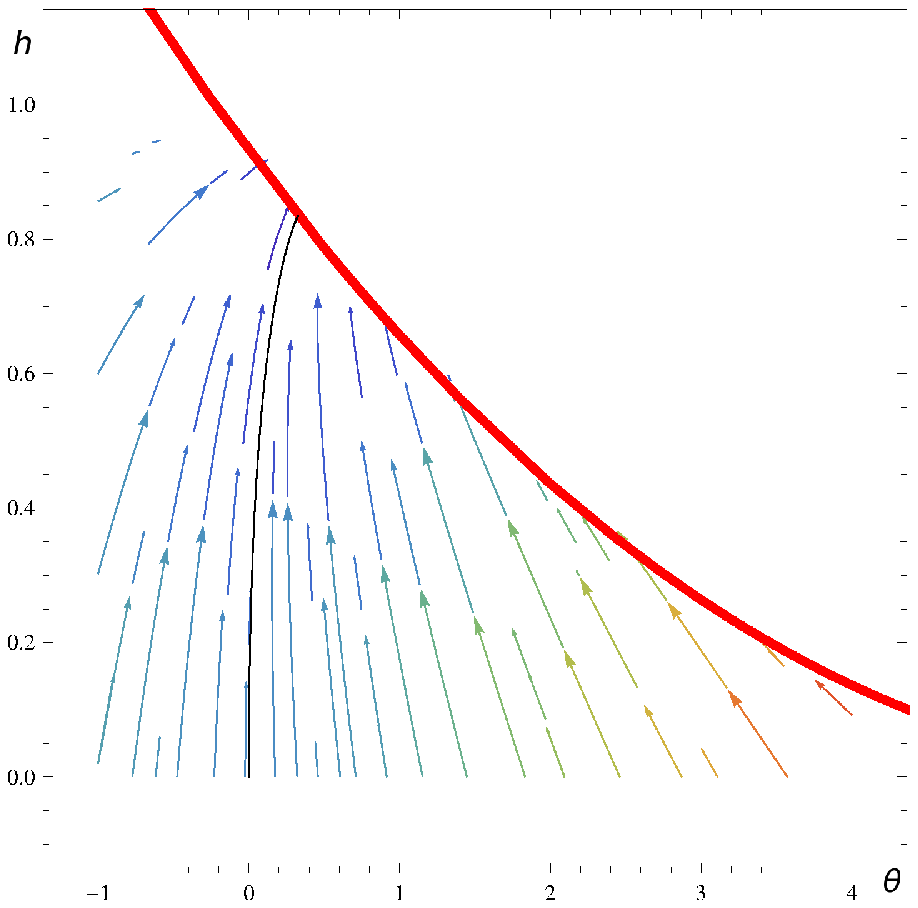}
		\caption{}\label{fig:SMdynamicsSnap}
	\end{subfigure}
	\caption{(a) The slow manifold (SM) and the break curve (thick line) of the system \eqref{eq:Atri3Da-Asym}--\eqref{eq:Atri3Dc-Asym}, plotted using, respectively, expressions \eqref{eq:SM} and \eqref{eq:FoldA}--\eqref{eq:FoldB}, when $T(c)=10c/(1+10c)$, $\mu=0.3$, and $\lambda=1$. (b) Phase plane of slow dynamics on the stable part of the SM ($\theta$-$h$ projection); the trajectories (arrows) stop at the break curve (thick line). (Trajectories would take off at the break curve to enter a fast motion.)}\label{fig:SMandFOLD}
\end{figure}

To determine the slow manifold (SM) we set $\epsilon=0$ in \eqref{eq:Atri3Da-Asym}, and obtain the surface
\begin{align}
F_1=0 \Rightarrow &\mu K_1h\frac{b+c}{1+c}-\frac{\Gamma c}{K+c}+\lambda \theta=0.
\nonumber\\
\label{eq:SM}
\Rightarrow &(\mu hK_1-\Gamma+\lambda \theta)c^2+(\mu hK_1(b+K)-\Gamma+\lambda \theta(K+1))c+(\mu hK_1bK+\lambda \theta)=0
\end{align}
\eqref{eq:SM} is a quadratic in $c$ with the roots $c_{-}(h,\theta)$ and $c_{+}(h,\theta)$, which correspond, respectively, to the stable and the unstable branch of the SM. These branches are joined by the ``break" curve, where the SM is tangent to the fast foliation, and which is defined by the conditions $\partial F_1/\partial c=0, F_1=0$. These conditions lead to the following parametric expressions for the break curve\footnote{Frequently, the break curve is called the fold curve but we do not use the latter term here in order to avoid confusion with the use of the term in Section \ref{sec:LinStab}}:
\begin{align}
\label{eq:FoldA}
h_F(c)&=\frac{\Gamma K}{\mu K_1}\frac{1}{1-b}\frac{(1+c)^2}{(K+c)^2},\\
\label{eq:FoldB}
\theta_F(c)&=\frac{\Gamma}{\lambda}\left(\frac{c}{K+c}-\frac{K}{1-b}\frac{(1+c)(b+c)}{(K+c)^2}\right).
\end{align}
For the stable part, 
$\partial F_1/\partial c<0$ and for the unstable part $\partial F_1/\partial c>0$.  The SM is also the set of equilibria of the fast system, and the stable and unstable parts of the SM correspond, respectively, to the stable and unstable steady states of the fast system. In Figure \ref{fig:SMandFOLD-BetterView}	 we plot the SM and the break curve (thick line) for $\mu=0.3$ and $\lambda=1$; for $\theta=0$ the SM is nullcline \eqref{eq:Nullcline1} of the Atri model. There is a slow motion on the stable part of the SM until the break curve is reached. The slow dynamics on the stable part of the SM are obtained by solving equations \eqref{eq:Atri3Db-Asym} and \eqref{eq:Atri3Dc-Asym} subject to the constraint $c=c_{-}(h,\theta)$. In Figure \ref{fig:SMdynamicsSnap} we plot the phase plane on the stable part of the SM. Trajectories stop at the break curve (plotted with a thick line). At the break curve, the trajectories take off into a fast motion---more details on this part of the motion are given below. To obtain the trajectories we can solve \eqref{eq:Atri3Db-Asym} and \eqref{eq:Atri3Dc-Asym}, an ODE system for the variables $\theta$ and $h$ subject to the constraint $F_1=0$\footnote{Alternatively, we can convert \eqref{eq:Atri3Db-Asym} and \eqref{eq:Atri3Dc-Asym} to an ODE system for the variables $c$ and $h$ by using $F_1=0$ to express $\theta$ in terms of $c$ and $h$, then differentiating with respect to $t$ to obtain $\theta_t$ in terms of  $c_t$ and $h_t$, and substituting back in \eqref{eq:Atri3Db-Asym}.}. 
$\partial F_1/\partial c<0$ and for the unstable part $\partial F_1/\partial c>0$.  The SM is also the set of equilibria of the fast system, and the stable and unstable parts of the SM correspond, respectively, to the stable and unstable steady states of the fast system. In Figure \ref{fig:SMandFOLD-BetterView}	 we plot the SM and the break curve (thick line) for $\mu=0.3$ and $\lambda=1$; for $\theta=0$ the SM is nullcline \eqref{eq:Nullcline1} of the Atri model. There is a slow motion on the stable part of the SM until the break curve is reached. The slow dynamics on the stable part of the SM are obtained by solving equations \eqref{eq:Atri3Db-Asym} and \eqref{eq:Atri3Dc-Asym} subject to the constraint $c=c_{-}(h,\theta)$. In Figure \ref{fig:SMdynamicsSnap} we plot the phase plane on the stable part of the SM. Trajectories stop at the break curve (plotted with a thick line). At the break curve, the trajectories take off into a fast motion until they hit the SM again. To obtain the trajectories we can solve \eqref{eq:Atri3Db-Asym} and \eqref{eq:Atri3Dc-Asym}, an ODE system for the variables $\theta$ and $h$ subject to the constraint $F_1=0$\footnote{Alternatively, we can convert \eqref{eq:Atri3Db-Asym} and \eqref{eq:Atri3Dc-Asym} to an ODE system for the variables $c$ and $h$ by using $F_1=0$ to express $\theta$ in terms of $c$ and $h$, then differentiating with respect to $t$ to obtain $\theta_t$ in terms of  $c_t$ and $h_t$, and substituting back in \eqref{eq:Atri3Db-Asym}.}. 

\subsubsection{$c$ large region -- `transition' layer}
When a trajectory leaves the stable part of the SM at the break curve, since there is no other stable branch on the SM to go to, $c$ becomes very large, and under a certain condition, the trajectory eventually crosses the SM again, and then moves towards the stable branch.  
To investigate this region we follow an asymptotic analysis similar to the one used in the Atri model previously, setting $c=\hat c/\epsilon$ and inserting the expansions $\hat c=\hat c_0+O(\epsilon)$, $\hat \theta=\hat \theta_0+O(\epsilon)$, $\hat  h=\hat  h_0+O(\epsilon)$, into the original system \eqref{eq:Atri3Da-Asym}--\eqref{eq:Atri3Dc-Asym}. We then obtain
\begin{align}
\label{eq:Atri3Da-III}
\frac{d\hat c_0}{dt}&=\mu \hat h_0 - \frac{\Gamma}{K_1}+\Lambda \hat \theta_0,\\
\label{eq:Atri3Db-III}
\frac{d \hat \theta_0}{dt}&=-\hat \theta_0+T_s\\
\label{eq:Atri3Dc-III}
\frac{d\hat  h_0}{dt}&=-\hat h_0,
\end{align}
assuming $T \to T_s$ as $c \to \infty$.
The latter reduced system can be solved analytically. From \eqref{eq:Atri3Db-III} and  \eqref{eq:Atri3Dc-III}  we obtain
\begin{align}
\label{eq:Atri3Dbc-III-Sol}
\hat \theta_0(t)=T_s+A_1e^{-t},\,\,\,\hat  h_0=\hat  h_{00}e^{-t},
\end{align}
where $A_1$ and $\hat h_{00}$ are constants of integration.
Matching as $t \to 0$, $\hat  h_0=\hat  h_{00}\to h_F$ and $\hat  \theta_0 \to \theta_F$, where $(\theta_F, h_F)$ is the point on the break curve where the trajectory is leaving the SM; as we have seen, this point is  parametrised by equations \eqref{eq:FoldA}  and \eqref{eq:FoldB}. Substituting \eqref{eq:Atri3Dbc-III-Sol} in \eqref{eq:Atri3Da-III} we find
\begin{align}
\label{eq:Atri3Da-III-Sol}
\frac{d\hat c_0}{dt}=(\mu h_F-\Lambda T_s)e^{-t}- \frac{\Gamma}{K_1}+\Lambda T_s+\Lambda \theta_F.
\end{align}
Integrating and setting the matching condition $\hat  c_0 \to 0$ as $t \to 0$, we find
\begin{align}
\label{eq:Atri3Da-III-Sol-Int}
&\hat c_0(t)=(\mu h_F-\Lambda T_s)(1-e^{-t})+(- \frac{\Gamma}{K_1}+\Lambda T_s +\Lambda \theta_F)t,\\
\label{eq:Atri3Db-III-Sol-Int}
&\hat h_0(t)=h_Fe^{-t},\,\,\hat \theta_0(t)=T_s(1-e^{-t})+\theta_Fe^{-t}.
\end{align}
In \eqref{eq:Atri3Da-III-Sol}, to determine the time needed for the trajectories to hit the SM for large $c$ we set $\frac{d\hat c_0}{dt}=0$. Solving for $t$ we find
\begin{align}
\label{eq:Atri3D-IIISolb}
t_{\rm  TURNING}=\ln\left(\frac{\mu K_1h_F-\lambda T_s}{\Gamma-\lambda T_s-\lambda \theta_F}\right).
\end{align}
Since $h_F$ scales with $1/(\mu K_1)$, $t_{\rm TURNING}$ does not depend on $\mu$ but, due to the presence of $\lambda T_s$ in \eqref{eq:Atri3D-IIISolb}, it does depend on the parameter $\lambda$. Comparing with the expression \eqref{eq:AtriND1IIISolb} for $t_{\rm TURNING}$ in the Atri model we conclude that the mechanical terms \emph{do affect} the dynamics at large $c$, at leading order. In particular, we find that in the mechanochemical model the trajectories turn only if 
\begin{align}
\label{eq:tTURN3D}
t_{\rm  TURNING}>0 \implies \frac{\mu K_1h_F-\lambda T_s}{\Gamma-\lambda T_s-\lambda \theta_F}>1 \implies \lambda \theta_F>\Gamma-\mu K_1h_F.
\end{align}
Substituting expressions  \eqref{eq:FoldA}--\eqref{eq:FoldB} in \eqref{eq:tTURN3D}, and performing the calculations, everything cancels out and \eqref{eq:tTURN3D} reduces to
\begin{align}K<1.
\end{align} 
The latter condition is satisfied for the chosen parameter values in this paper since $K=1/7$, but for $K \ge 1$ the model predicts that some trajectories would, unphysically, shoot off to infinity after leaving the slow manifold at the break curve.
 
Note, finally, that from \eqref{eq:Atri3Da-III} we also have $\frac{d\hat c_0}{dt}=0 \Rightarrow \mu \hat h_0 = \frac{\Gamma}{K_1}-\Lambda \hat \theta_0$, which gives the plane in the $c,h,\theta$ space where the trajectories cross the SM, depending on where they start off when they leave the break curve. 

\section{Summary, conclusions and discussion}
\label{sec:Conclusions}
In this paper we have first revisited the well-known Atri model of calcium dynamics \cite{ atri1993single} and analysed its bifurcation structure in detail as the bifurcation parameter $\mu$ (a function of the $IP_3$ concentration) was increased; we identified some new bifurcations that are very sensitively dependent on $\mu$. We then analysed the Atri model asymptotically for the range of $\mu$ for which the system exhibits relaxation oscillations, assuming that $c$ was a fast variable compared to $h$, the percentage of active $IP_3$ receptors. This assumption was supported by the parameter values given in \cite{ atri1993single}. The asymptotic analysis has qualitative differences from the asymptotic analysis of Van der Pol-like systems since nullcline \eqref{eq:Nullcline1} of the Atri system is not a cubic but asymptotes to a constant value for $c$ large; in particular we identified a ``transition layer"' for $c$ large, where the trajectory crosses the nullcline on the unstable branch and instantaneously reverses direction --this corresponds to a very 'spiky' oscillation. 

Additionally, inspired by experimental results which show that cells release calcium in response to mechanical stimuli but also that calcium release causes the cells to contract/expand we have developed a new, simple mechanochemical, three-dimensional ODE model that consists of an ODE for $\theta$, the cell/tissue dilatation, derived consistently from a full viscoelastic \emph{ansatz}, and the two Atri ODEs. An important feature of this model is a \emph{two-way} mechanochemical feedback mechanism which, to our knowledge, has not been employed previously in the calcium literature. In the viscoelastic ansatz we have used, we assumed that the cell/tissue is a linear, one-dimensional, Kelvin-Voigt viscoelastic material and that the stress tensor contains a calcium-dependent term, modelled by the function $T(c)$. We also modified the Atri calcium ODE by adding a new ``stretch-activation", mechanical source term $\lambda \theta$; $\lambda$ was treated as a second bifurcation parameter in the analysis. We then linearised the ODE system and developed a parametric method valid  for \emph{any} functional form of $T(c)$ with which we can easily plot the key bifurcation curves of the system in the $\mu$--$\lambda$ plane. We two Hill functions, $T(c)=\alpha c/(1+\alpha c)$ and $T(c)=\alpha c^2/(1+\alpha c^2)$, we plotted the Hopf curves since their interior is the $\mu$-$\lambda$ parameter region for which relaxation oscillations are sustained. For both these cases of $T$ as $\lambda$ was increased the relaxation oscillations were eventually suppressed at a value of $\lambda$, $\lambda_{\rm max}$. This was a key result of this paper, and it would be interesting to devise an appropriate experiment to investigate its validity. Note that as $\alpha$, the rate of growth of $T$ at $c=0$, was decreased $\lambda_{\rm max}$ increases for both cases of $T(c)$. 

Furthermore, choosing $T(c)=10 c/(1+10 c)$ we used XPPAUT to study the amplitude and frequency of the oscillations as $\lambda$ and $\mu$ were increased. We found that: (i) for relatively low values of $\lambda$ the amplitude increases rapidly with $\mu$ near the two Hopf points but slowly in an intermediate region (see Figures \ref{fig:AtriMechsLam0Pt1}--\ref{fig:AtriMechsLam1}). The intermediate region decreases with $\lambda$ as the LHP and the RHP approach each other; at a sufficiently large value of $\lambda$ RHP becomes supercritical (see Figure \ref{fig:AtriMechsLam1Pt5}) and the amplitude changes slowly for all values of $\mu$ (ii) For a fixed value of $\mu$, e.g. $\mu=0.3$ we see in Figures \ref{fig:AtriMechsLam0Pt1} and \ref{fig:AtriMechsLam0Pt5} that the amplitude stays constant with $\lambda$ (iii) the oscillation frequency rapidly increases close to the LHP and the RHP and there is an intermediate region where the frequency varies slowly with $\mu$ for smaller values of $\lambda$ (see Figures \ref{fig:Lam0Pt1Freq}--\ref{fig:Lam1Freq}). This behaviour is qualitatively similar to the $\lambda=0$ case (Atri model, see Figure \ref{fig:AtriBifnFrequency}). The behaviour changes when $\lambda$ further increases though the intermediate region seems to have vanish (see Figure \ref{fig:Lam1Pt5Freq}).

Finally, we have developed an asymptotic analysis of the mechanochemical model when calcium dynamics are fast compared to the dynamics of the other two variables. We assumed $T(c) \to T_s$, a constant, as  $c$ gets large. We analysed the slow manifold and described a typical periodic orbit trajectory, at leading order, as a combination of three lower-dimensional ODE systems which can be solved very easily; a 2D slow (nonlinear) system, a 1D fast (nonlinear) system, and a 2D linear system in the ``transition region'' where $c$ is large. We, thus, obtained insights about the behaviour of our 3D mechanochemical model by studying very simple, lower-dimensional systems. In particular, the dynamics in the transition region were solved for analytically; and identifying an expression for $t_{\rm TURNING}$ (see \eqref{eq:Atri3D-IIISolb}) we established that, for the parameter values of \cite{atri1993single}, and any value of $\lambda$ there are no `unphysical' trajectories that `escape' to infinity. Previous mechanochemical models \cite{warren2010, kobayashi2014, yao2016} were higher-dimensional and could only be analysed numerically, which limited the understanding of the underlying mechanochemical feedback mechanism. Therefore, we have introduced a simple, mechanochemical model that can easily analyse with semi-analytical methods but which still captures important aspects of the mechanochemical processes identified in experiments. It is indeed our hope that the proposed model can inspire more experimental work.

\begin{acknowledgements}
We thank James Sneyd, Vivien Kirk, and Ruediger Thul for valuable discussions on calcium modelling, and Bard Ermentrout for promptly providing support with the XPPAUT continuation software. Katerina Kaouri also acknowledges support from an STSM grant awarded by COST Action TD1409 (Mathematics for Industry Network, MI-NET) for a research visit to Oxford University.
\end{acknowledgements}

\bibliographystyle{spmpsci}      
\bibliography{refsCW-Aug16}   

\end{document}